\documentclass{amsart}
\usepackage{amssymb}
\usepackage{amsfonts}
\usepackage{latexsym}

\newtheorem{theorem}{Theorem}[section]
\newtheorem{lemma}[theorem]{Lemma}
\newtheorem{proposition}[theorem]{Proposition}
\newtheorem{corollary}[theorem]{Corollary} 
\theoremstyle{definition}  
\newtheorem{definition}[theorem]{Definition}

\newtheorem{example}[theorem]{Example}

\newtheorem{remark}[theorem]{Remark}

\newcommand{\e}{{\rm e}}

\newcommand{\End}{\text{End}}

\newcommand{\h}{\mathfrak{h}}

\newcommand{\ben}{\begin{enumerate}}
\newcommand{\een}{\end{enumerate}}

\begin{document}
\title{Cherednik and Hecke algebras of varieties with a finite
group action}
\author{Pavel Etingof}
\address{Department of Mathematics, Massachusetts Institute of
Technology,
Cambridge, MA 02139, USA}
\email{etingof@math.mit.edu}

\maketitle

\section{Introduction}

This paper is an expanded and updated version of the unpublished 2004 preprint \cite{Et1}. 
It includes a more detailed description of the basics of the  theory of Cherednik and Hecke algebras of varieties started in \cite{Et1}, as well as a new Section 4, which reviews the developments in this theory since 2004 with references to the relevant literature.  

Let $\h$ be a finite dimensional complex vector space, 
and $G$ be a finite subgroup of $GL(\h)$. 
To this data one can attach a family of algebras $H_{t,c}(\h,G)$,
called the rational Cherednik algebras (see \cite{EG});
for $t=1$ it provides the universal deformation of $G\ltimes
D(\h)$ (where $D(\h)$ is the algebra of differential operators on $\h$). 
These algebras are generated by $G,\h,\h^*$ with 
certain commutation relations, and are parametrized 
by pairs $(t,c)$, where $t$ is a complex number, and $c$ is a
conjugation invariant function on the set of complex reflections in $G$. 
They have a rich representation theory  
and deep connections with combinatorics (Macdonald theory, $n!$
conjecture) and algebraic geometry (Hilbert schemes, resolutions  
of symplectic quotient singularities). 

The purpose of this paper is to introduce ``global'' 
analogues of rational Cherednik algebras, attached to 
any smooth complex algebraic variety $X$ with an action of a finite group $G$;
the usual rational Cherednik algebras are recovered 
in the case when $X$ is a vector space and $G$ acts linearly.
 
More specifically, let $G$ be a finite group of automorphisms 
of $X$, and let $S$ be the set of pairs $(Y,g)$, where $g\in G$, and $Y$
is a connected component of the set $X^g$ of $g$-fixed points in
$X$ which has codimension 1 in $X$ (we will call such a component a {\it reflection hypersurface}). 
Suppose that $X$ is affine. Then we define (in Section 2 of the paper) 
a family of algebras $H_{t,c,\omega}(X,G)$,
where $t\in \Bbb C$, $c$ is a $G$-invariant function on $S$, and $\omega$ is a
$G$-invariant closed 2-form on $X$. This family for
$t=1$ provides a universal deformation of the algebra
$H_{1,0,0}(X,G)=G\ltimes D(X)$, where $D(X)$ is the algebra of
differential operators on $X$ (assuming that $\omega$ runs
through a space of forms bijectively representing $H^2(X,\Bbb C)^G$).  

If $X$ is not affine, then we define a sheaf of algebras 
$H_{t,c,\omega,X,G}$ rather
than a single algebra. In this case the parameters are the same,
except that $\omega$ runs over a space representing classes 
of $G$-equivariant twisted differential operator (tdo) algebras
on $X$ (see \cite{BB}, Section 2). 

We find that much of the theory of rational Cherednik algebras
survives in the global case. In particular, one can
define the spherical subalgebra, which is both commutative 
and isomorphic to the center of the Cherednik algebra
in the case $t=0$. The spectrum of this algebra
is ``the Calogero-Moser space'' of $X$, which is a global version
of the similar space defined in \cite{EG}. This includes, in
particular, Calogero-Moser spaces attached to 
symmetric powers of algebraic curves. One can also define 
the global analog of the theory of quasiinvariants 
which was worked out in \cite{FV,EG,BEG} and references therein.
These results can be generalized to the case when $X$ is 
a complex analytic variety. 

In Section 3, we discuss an application of the theory of global 
Cherednik algebras for complex analytic varieties. 
Namely, for an analytic $G$-variety $X$ we define its Hecke algebra,
which is a certain explicitly defined 
formal deformation of the group algebra of the orbifold 
fundamental group of $X/G$. We show that if \linebreak $\pi_2(X)\otimes
\Bbb Q=0$ (a condition that cannot be removed), then the Hecke
algebra is a flat deformation. This includes usual, affine, and
double affine Hecke algebras for Weyl groups, as well as Hecke
algebras for complex reflection groups. The proof is based 
on showing that the regular representation 
of the orbifold fundamental group can be deformed to a
representation of the Hecke algebra. The required deformation is
constructed by applying the KZ functor to a module over 
the global Cherednik algebra. 

In Section 4 we review the developments in the theory of Cherednik and Hecke algebras of varieties with a finite group action since their introduction in the preprint \cite{Et1}, and give the corresponding references. This section shows how these algebras fit into a bigger representation-theoretic, deformation-theoretic, and geometric picture.  

{\bf Acknowledgements.} This paper is dedicated to the the 80-th anniversary of V. I. Arnold. His books, courses, and his approach to mathematics in general have profoundly influenced me since my youth and will continue to do so in the future. 

I am very grateful to A. Braverman for many explanations
about D-modules, in particular for proofs of 
Propositions \ref{groth} and \ref{hoh}. I thank V. Ginzburg 
and A. Okounkov for useful discussions, and D. Thompson for comments on the draft of this paper. 
My work was partially supported by the NSF grants DMS-9988796, DMS-1502244, and the CRDF grant RM1-2545-MO-03.

\section{The Cherednik algebra of a $G$-variety}

\subsection{Basic results about $D$-modules}

In this subsection we will prove a few basic results about D-modules, to be used below. We refer the reader to the books \cite{Bo,HTT} for basics on D-modules. 

Throughout the paper, we work over $\Bbb C$. Let $X$ be a smooth algebraic variety. Let ${\mathcal{O}}_X$ be the structure sheaf of $X$, and $D_X$ be the sheaf of differential operators on $X$. We will need the following
well known result. 

\begin{proposition} \label{groth} One has a natural isomorphism ${\rm Ext}^i_{D_X}({\mathcal
O}_X,{\mathcal O}_X)\cong H^i(X,\Bbb C)$. 
\end{proposition}

\begin{proof} Let $D^\bullet_X$ be the dual D-module De Rham
complex of $X$: 
$$
D_X\otimes_{{\mathcal O}_X}\Omega^{-d}_X\to...\to
D_X\otimes_{{\mathcal O}_X}\Omega^{-1}_X\to
D_X,
$$
where $\Omega^{-i}_X$ is the sheaf of polyvector fields 
of rank $i$ on $X$. Then $D^\bullet_X$ is a locally projective 
resolution of the D-module ${\mathcal O}_X$. 
Thus, the required Ext groups are equal to the hypercohomology
groups of the complex of sheaves ${\rm Hom}_{D_X}
(D^\bullet_X,{\mathcal O}_X)$. This complex is just the 
algebraic De Rham complex of $X$. But by Grothendieck's algebraic De Rham theorem, 
the hypercohomology of the De Rham complex is naturally isomorphic to the
cohomology of $X$. This proves the proposition.
\end{proof} 

\begin{remark}\label{analvar} Proposition \ref{groth} also holds for analytic 
varieties, with the same proof (using the usual De Rham theorem for complex manifolds). 
\end{remark} 

Now let $X$ be a smooth algebraic variety over $\Bbb C$, and 
$g$ be an automorphism of $X$ of finite order. 
In this case the set of fixed points $X^g$ is a smooth algebraic
variety, consisting of finitely many (say, $N_g$) connected components 
$X_j^g$ (possibly of different dimensions).  

Assume that $X$ is affine, and let $D(X)$ be the algebra of
differential operators on $X$. Then $g$ acts naturally as an automorphism 
of $D(X)$; namely, this action is induced by the action of $g$ on ${\mathcal O}(X)$ given by $(gf)(x):=f(g^{-1}x)$. Thus we can define the bimodule $D(X)g$ over $D(X)$
in an obvious way: $a\circ bg\circ c:=abg(c)g$ for $a,b,c\in D(X)$.  

The following proposition computes the Hochschild cohomology of $D(X)$ with
coefficients in $D(X)g$. 

\begin{proposition}\label{hoh}
We have a natural isomorphism 
$$
HH^m(D(X),D(X)g)\cong\oplus_{j=1}^{N_g}H^{m-2{\rm
codim}X_j^g}(X_j^g,\Bbb C).
$$
In particular, this isomorphism is equivariant with respect to
any automorphism of $X$ that commutes with $g$. 
\end{proposition}

\begin{proof}
Let $\widetilde{D}_X$ and $\widetilde{D}_Xg$ 
be the left D-modules on $X\times X$ corresponding to $D(X)$ and $D(X)g$, respectively, under the standard equivalence 
between left and right D-modules on the second factor $X$, see
e.g. \cite{Bo,HTT}. Let $\Delta: X\to X\times X$
be the diagonal map, and $i: X\to X\times X$ be given by the
formula $x\to (x,g^{-1}x)$. Then $\widetilde{D}_X=\Delta_*\mathcal{O}_X$ and 
$\widetilde{D}_Xg=i_*\mathcal{O}_X$. Thus, using that $\Delta_*$ is left adjoint of $\Delta^!$, we have 
$$
HH^m(D(X),D(X)g)={\rm Ext}^m_{D(X)\otimes D(X)^{\rm op}}(D(X),D(X)g)=
{\rm Ext}^m_{D_{X\times X}}(\widetilde{D}_X,\widetilde{D}_Xg)=
$$
$$
=
{\rm Ext}_{D_{X\times X}}^m
(\Delta_*{\mathcal O}_X,i_*{\mathcal O}_X)=
{\rm Ext}_{D_X}^m({\mathcal O}_X,\Delta^!i_*{\mathcal O}_X).
$$

Let $\eta_j: X^g_j\to X$ be the tautological
embedding. Applying base change to the composition $\Delta^!i_*$ and using that $\eta_j^!\mathcal{O}_X=O_{X_j^g}[-{\rm codim}X_j^g]$, we have 
$$
{\rm Ext}_{D_X}^m
({\mathcal O}_X,\Delta^!i_*{\mathcal O}_X)=
{\rm Ext}_{D_X}^m
({\mathcal O}_X,\oplus_{j=1}^{N_g}\eta_{j*}{\mathcal
O}_{X_j^g}[-{\rm codim}X_j^g]).
$$
Since $\eta_j^*$ is left adjoint to $\eta_{j*}$, and 
$\eta_j^*\mathcal{O}_X=\mathcal{O}_{X_j^g}[+{\rm codim} X_j^g]$, we have  
$$
{\rm Ext}_{D_X}^m
({\mathcal O}_X,\oplus_{j=1}^{N_g}\eta_{j*}{\mathcal
O}_{X_j^g}[-{\rm codim}X_j^g])=
$$
$$
\oplus_{j=1}^{N_g}{\rm Ext}_{D_{X_j^g}}^m
(\eta_j^*{\mathcal O}_X,{\mathcal
O}_{X_j^g}[-{\rm codim}X_j^g])=
$$
$$
\oplus_{j=1}^{N_g}{\rm Ext}_{D_{X_j^g}}^m
({\mathcal O}_{X_j^g}[+{\rm codim}X_j^g],
{\mathcal O}_{X_j^g}[-{\rm codim}X_j^g])=
$$
$$
=\oplus_{j=1}^{N_g}{\rm Ext}_{D_{X_j^g}}^{m-2{\rm codim}X_j^g}
({\mathcal O}_{X_j^g},{\mathcal O}_{X_j^g})=
\oplus_{j=1}^{N_g}H^{m-2{\rm
codim}X_j^g}(X_j^g,\Bbb C),
$$
where the last equality follows from Proposition
\ref{groth}. This implies the required statement. 
\end{proof}

For $X=\Bbb A^n$, Proposition \ref{hoh} appears in \cite{AFLS}. 
We also note that an analog of Proposition \ref{hoh} for $g=1$ and smooth real manifolds 
is due to Kassel and Mitschi (see e.g. \cite{BG}). 

\begin{corollary} \label{cohGD}
We have a natural isomorphism
$$
HH^m(G\ltimes D(X),G\ltimes D(X))\cong
(\oplus_{g\in G}\oplus_{j=1}^{N_g}H^{m-2{\rm
codim}X_j^g}(X_j^g,\Bbb C))^G.
$$
\end{corollary}

\begin{proof}
It is well known (see e.g. \cite{AFLS}) that for any 
algebra $A$ over a field of characteristic zero
with an action of a finite group $G$, one has
$$
HH^m(G\ltimes A,G\ltimes A)=(\oplus_{g\in G}HH^m(A,Ag))^G.
$$
Indeed, 
$$
HH^m(G\ltimes A,G\ltimes A)=
{\rm Ext}^m_{(G\times G)\ltimes A\otimes A^{op}}(G\ltimes
A,G\ltimes A)
$$
By Shapiro's lemma, this equals 
$$
{\rm Ext}^m_{G_{\rm diagonal}\ltimes A\otimes A^{op}}(A,G\ltimes
A),
$$
which (by Maschke's theorem for $G$) is equal to
$$
{\rm Ext}^m_{A\otimes A^{op}}(A,G\ltimes A)^G=
HH^m(A,G\ltimes A)^G,
$$
as desired. 

So for $A=D(X)$, we find
$$
HH^m(G\ltimes D(X),G\ltimes D(X))=
(\oplus_{g\in G} HH^m(D(X),D(X)g))^G.
$$
The rest follows from Proposition \ref{hoh}.
\end{proof} 

\subsection{Twisted differential operators}\label{tdo}

Let us recall the theory of twisted differential
operators (see \cite{BB}, section 2).

Let $X$ be a smooth affine algebraic variety over $\Bbb C$.  
Let $D(X)$ be the algebra of algebraic differential operators on $X$.
Given a closed 2-form $\omega$ on $X$, we can define a
two-cocycle on the Lie algebra ${\rm Vect}(X)$ 
with coefficients in the module of regular functions ${\mathcal O}(X)$, 
given by $v,w\to \omega(v,w)$. This cocycle defines an abelian
extension of ${\rm Vect}(X)$ by ${\mathcal O}(X)$, which has an obvious 
structure of a Lie algebroid ${\Bbb L}_\omega$ over $X$. 
It is clear that this Lie algebroid depends only on 
the cohomology class of $\omega$, up to an isomorphism. 

Let $U=U({\Bbb L}_\omega)$ be the universal enveloping algebra 
of this Lie algebroid. For any $f\in {\mathcal O}(X)$, let $\widehat f$ 
be the corresponding section of ${\Bbb L}_\omega$. 
Let $I$ be the ideal in $U$ generated by the elements $\widehat
f-f$. The quotient $U/I$ is the algebra of twisted differential
operators $D_\omega(X)$. 

More explicitly, $D_\omega(X)$ can be defined 
as the algebra generated by ${\mathcal O}(X)$ and 
``Lie derivatives'' ${\bold L}_v$, $v\in {\rm Vect}(X)$,
with defining relations 
$$
f{\bold L}_v={\bold L}_{fv},\ [{\bold L}_v,f]=L_vf,\ [{\bold L}_v,{\bold L}_w]={\bold L}_{[v,w]}+\omega(v,w).
$$
This algebra depends only on the
cohomology class $[\omega]$ of 
$\omega$ (up to an isomorphism), and equals $D(X)$ if $\omega=0$. 

An important special case of twisted differential operators is 
the algebra of differential operators on a line bundle. 
Namely, let $L$ be a line bundle on $X$. Since $X$ is affine, $L$
admits an algebraic connection $\nabla$ with curvature $\omega$,
which is a closed 2-form on $X$. Then it is easy to show that 
the algebra $D(X,L)$ of differential operators on $L$ 
is isomorphic to $D_\omega(X)$. 

The classical analogs of twisted differential operator algebras
are twisted cotangent bundles. In the case of an affine variety
$X$, a twisted cotangent bundle is the 
usual cotangent bundle variety $T^*X$ equipped
with the new symplectic structure $\Omega'=\Omega+\pi^*\omega$, where 
$\Omega$ is the usual symplectic structure on $T^*X$
and $\pi: T^*X\to X$ is the projection. Here $\omega$ is a closed
2-form on $X$. We will denote the twisted cotangent bundle by
$T^*_\omega X$. 

If the variety $X$ is smooth but not necessarily affine, 
then sheaves of algebras of twisted differential operators 
$D_{\psi,X}$ and twisted cotangent bundles $T^*_\psi X$ are classified by 
elements $\psi$ of the hypercohomology $H^2(X,\Omega_X^{\ge 1})$, where 
$\Omega_X^{\ge 1}$ is the two-step complex 
of sheaves $\Omega^1_X\to \Omega^{2,cl}_X$, given by the De Rham
differential acting from 1-forms to closed 2-forms (sitting in
degrees 1 and 2, respectively). 
If $X$ is projective then this space is isomorphic 
to $H^{2,0}(X,\Bbb C)\oplus H^{1,1}(X,\Bbb C)$.
We note that if $X$ is not affine then 
the twisted cotangent bundle is in general not a vector bundle
but an affine space bundle, which need not be isomorphic to
$T^*X$. We refer the reader to \cite{BB},
Section 2, for details. 

\subsection{$D_\omega(X)$ as a universal deformation.} 

Let $E$ be a subspace of the space of closed 2-forms
on $X$ which projects isomorphically to $H^2(X)$.\footnote{Unless otherwise specified, cohomology of 
varieties is with complex coefficients.}  
Then $D_\omega(X)$, $\omega\in E$, is a family of algebras parametrized 
by $H^2(X)$. 

\begin{lemma}\label{nont}
Let $\omega\in E$. 
If the first order deformation 
$D_{\hbar\omega}(X)$ of $D(X)$ 
over $\Bbb C[\hbar]/\hbar^2$ is trivial, then $\omega=0$. 
\end{lemma}

\begin{proof}
If the deformation is trivial then the module ${\mathcal O}(X)$ 
over $D(X)$ can be lifted to this deformation. 
This lifting must be trivial as a
deformation of ${\mathcal O}(X)$-modules, i.e. it is isomorphic 
to ${\mathcal O}(X)[\hbar]/\hbar^2$ with the usual action of
functions, and the action of vector fields by 
$$
{\bold L}_vf=L_vf+\hbar A_v(f)
$$
for some operator $A_v$, where $L_vf$ is the usual Lie
derivative. Since $[{\bold L}_v,f]=L_vf$ in $D_{\hbar\omega}(X)$, we 
find that $A_v(f)$ commutes with operators of multiplication by
functions, so using that ${\bold L}_{gv}=g{\bold L}_v$, we get that
$A_v(f)=\eta(v)f$, where $\eta$ is a 1-form on $X$. 
From this, using the commutation relation between 
${\bold L}_v$ and ${\bold L}_w$, by a simple calculation we obtain
$d\eta=\omega$. Since the map $E\to H^2(X)$ is an isomorphism,
$\omega=0$. The lemma is proved.
\end{proof}

By Proposition \ref{hoh}, $HH^2(D(X),D(X))=H^2(X)$. Thus
using Lemma \ref{nont}, we see that the 1-parameter formal deformations 
of $D(X)$ induced by $D_\omega(X)$, $\omega\in E$,
 represent all elements of $H^2(X)$. 
Therefore, we get the following (apparently, well known)
theorem. 

\begin{theorem}\label{tdot} $D_\omega(X)$, where $\omega$ lies in the formal
neighborhood of the origin in $E$, is a universal formal
deformation of $D(X)$. 
\end{theorem} 

\begin{remark} In fact, it is easy to show that the natural map $E\to HH^2(D(X),D(X))\cong H^2(X)$ 
induced by the deformation $D_\omega(X)$ is given by the formula $\omega\mapsto [\omega]$. 
\end{remark} 

\subsection{Algebro-geometric preliminaries}\label{agp}
Let $Z$ be a smooth hypersurface in a smooth variety $X$. 
Let $i: Z\hookrightarrow X$ be the corresponding closed embedding. 
Let $N$ denote the normal bundle of $Z$ in $X$
(a line bundle). Let ${\mathcal O}_X(Z)$ 
denote the sheaf of regular functions on
$X\setminus Z$ which have a pole of at most first order at $Z$.
Then we have a natural map of ${\mathcal O}_X$-modules 
$\phi: {\mathcal O}_X(Z)\to i_*N$. Indeed, 
we have a natural residue map 
$\eta: {\mathcal O}_X(Z)\otimes_{{\mathcal O}_X}\Omega^1_X\to 
i_*{\mathcal O}_Z$ (where $\Omega^1_X$ is the sheaf of 1-forms on $X$, and $i_*$ denotes the direct image of 
quasicoherent sheaves), hence a map $\eta': {\mathcal O}_X(Z)\to 
i_*{\mathcal O}_Z\otimes_{{\mathcal O}_X}TX=i_*(TX|_Z)$ (where $TX=\Omega_X^{-1}$ is
the tangent bundle). The map $\phi$ is obtained by composing
$\eta'$ with the natural projection $TX|_Z\to N$.  

We have an exact sequence 
of ${\mathcal O}_X$-modules: 
$$
0\to {\mathcal O}_X\to {\mathcal O}_X(Z)\to i_*N\to 0
$$
(the third map is $\phi$). Thus we have a 
natural surjective map of ${\mathcal O}_X$-modules 
$\xi_Z: TX\to {\mathcal O}_X(Z)/{\mathcal O}_X$.

\subsection{Cartan's lemma} 

We will need the following well known lemma, due to H. Cartan. 

\begin{lemma}\label{cartlem}
(i) Let $X$ be a complex manifold with a holomorphic action of a finite group $G$, and $x\in X$ be a fixed point of $G$. 
Then there exists $G$-invariant open subsets $x\in U\subset X$, 
$0\in U_0\subset T_xX$ and a $G$-equivariant isomorphism $f: U_0\cong U$.
In other words, locally near $x$ the $G$-action can be linearized. 

(ii) (formal Cartan's lemma) Any action of a finite group $G$ on a formal polydisk over a field of characteristic zero 
is equivalent to a linear action.  
\end{lemma} 

\begin{proof} For (i) see \cite{C}, p.97. For (ii), see e.g. \cite{EM1}, Lemma 7.8. 
\end{proof}

\subsection{Rational Cherednik algebras}\label{rca}

Recall the definition of rational Cherednik algebras (see e.g. 
\cite{EG}). 

Let $\h$ be a finite dimensional complex vector space, 
and $G$ is a finite subgroup of $GL(\h)$. Let $S$ be the set of
complex reflections in $G$, i.e., elements which have only one
eigenvalue not equal to $1$. Let $t\in \Bbb C$, and $c: S\to \Bbb C$ be 
a $G$-invariant function. To this data one attaches an algebra
$H_{t,c}(\h,G)$, called the rational Cherednik algebra 
of $\h,G$. It is generated by $G,\h,\h^*$ with 
defining relations
$$
gxg^{-1}={}^gx,\ gyg^{-1}={}^gy,\ 
[x,x']=[y,y']=0,
$$
and the main commutation relation
$$
[y,x]=t(y,x)-\sum_{s\in S}c_s(y,\alpha_s)(x,\alpha_s^\vee)s,
$$
where $x,x'\in\h^*, y,y'\in \h$, $\alpha_s$ is 
a nonzero linear function on $\h$ vanishing on 
the fixed hyperplane of $s$ in $\h$, and $\alpha_s^\vee$ 
is the element of $\h$ vanishing on the fixed hyperplane of $s$
in $\h^*$, such that $(\alpha_s,\alpha_s^\vee)=2$.

The algebra $H_{t,c}(\h,G)$ 
has a natural representation in $M=\Bbb C[\h]$ 
by Dunkl-Opdam operators \cite{DO}:
$$
x\mapsto x,\ g\mapsto g,\
y\mapsto D_y:=t\frac{\partial}{\partial y}
+\sum_{s\in S}\frac{2c_s}{1-\lambda_s}\frac{(\alpha_s,y)}{\alpha_s}(s-1)
$$
($g\in G,x\in \h^*,y\in \h$), where $\lambda_s$ 
is the nontrivial eigenvalue of $s$ in $\h^*$. 

This representation is faithful. Thus we 
can alternatively define $H_{t,c}(\h,G)$ (for $t,c$ being
variables) as follows: 
$H_{t,c}(\h,G)$ is the subalgebra 
of the algebra \linebreak $G\ltimes D(\h)_r[t,c]$ generated over
$\Bbb C[t,c]$ 
by $\Bbb C[\h]=S\h^*$, $G$, and the operators $D_y$
(here $D(\h)_r$ is the algebra of differential operators on $\h$ 
with rational coefficients). 

The algebra $H_{t,c}(\h,G)$ has an 
increasing filtration $F^\bullet$, defined by 
the rule $\deg(G)=\deg(\h^*)=0$, $\deg(\h)=1$. 
The PBW theorem (see e.g. \cite{EG} or \cite{EM1}, Subsection 3.2)
says that ${\rm gr}_F(H_{t,c}(\h,G))=G\ltimes
\Bbb C[\h\oplus \h^*]=G\ltimes S(\h^*\oplus \h)$. 

Let $\lbrace{y_i\rbrace}$ be a basis of $\h$ and $\lbrace{x_i\rbrace}$ be the dual basis of
$\h^*$. Define the Euler element 
$$
\bold h=\sum_i x_iy_i+\frac{\ell t}{2}-\sum_{s\in
S}\frac{2c_s}{1-\lambda_s}s\in H_{t,c}(\h,G),
$$
where $\ell=\dim \h$ (see \cite{GGOR}). 
Recall that $[\bold h,x]=x$ and $[\bold h,y]=-y$ for 
$x\in \h^*$, $y\in \h$.

Recall \cite{GGOR} that the category ${\mathcal O}_c(\h,G)$ 
for $H_{1,c}(\h,G)$ is the category of representations 
of $H_{1,c}(\h,G)$ which are direct sums 
of finite dimensional generalized eigenspaces of $\bold h$, 
with real part of the spectrum of $\bold h$ bounded from below. 
Thus, any module from ${\mathcal O}_c(\h,G)$ is graded by generalized
eigenvalues of $\bold h$. Note also that any $M\in {\mathcal O}_c(\h,G)$
is finitely generated over the subalgebra $\Bbb C[\h]$, and that 
${\mathcal O}_c(\h,G)$ may be alternatively defined as the category of 
$H_{1,c}(\h,G)$-modules which are finitely generated over $\Bbb C[\h]=S\h^*$ and 
have a locally nilpotent action of $\h$.  

\subsection{The formal completion of the rational Cherednik algebra}\label{comple} 

Now consider the degree-wise formal completion $\widehat H_{t,c}(\h,G)$
of the rational Cherednik algebra $H_{t,c}(\h,G)$, i.e., its restriction to the formal neighborhood of zero in
$\h$ as a $\Bbb C[\h]$-module: $\widehat H_{t,c}(\h,G):=\Bbb C[[\h]]\otimes_{\Bbb C[\h]}H_{t,c}(\h,G)$. This completion has a natural algebra structure, and comes with an increasing filtration $F^\bullet$ defined by the rule 
$F^i\widehat H_{t,c}(\h,G)=\Bbb C[[\h]]\otimes_{\Bbb
C[\h]}F^iH_{t,c}(\h,G)$. 

If $M\in {\mathcal O}_c(\h,G)$,
then the completion $\widehat M$ of $M$ by the grading defined by eigenvalues
of ${\bold h}$ is a representation of the algebra 
$\widehat H_{t,c}(\h,G)$. 
The category $\widehat{\mathcal O}_c(\h,G)$ 
 is defined as the category of modules over the algebra
$\widehat H_{t,c}(\h,G)$ of the form $\widehat M$.
Then the functor ${\mathcal O}_c(\h,G)\to \widehat{\mathcal O}_c(\h,G)$ given by $M\mapsto \widehat M$ is an equivalence of categories (see \cite{BE}, Theorem 2.3). 

We will need the following result. 

\begin{proposition}\label{o=fin} (see also \cite{BE}, Theorem 2.3)
The category $\widehat{\mathcal O}_c(\h,G)$ coincides with the
category of $\widehat H_{1,c}(\h,G)$-modules
which are finitely generated over $\Bbb C[[\h]]$. 
\end{proposition}

\begin{proof} It is clear that the first category is a full
subcategory of the second one; so our job is to show that any
 $\widehat H_{1,c}(\h,G)$-module $N$ which is 
finitely generated over $\Bbb C[[\h]]$ belongs to ${\mathcal O}_c(\h,G)$.

Let $I$ be the maximal ideal in $\Bbb C[[\h]]$. 
The module $N$ has a decreasing filtration $N\supset IN\supset
I^2N\supset...,$, such that the quotients $N/I^kN$ are finite
dimensional, and $N$ is the inverse limit of $N/I^kN$. 
The element ${\bold h}$ normalizes $I$, and ${\rm ad}{\bold h}$
has positive integer eigenvalues on $I$. This implies that 
for any $\lambda\in \Bbb C$ the generalized eigenspaces
$(N/I^kN)(\lambda)$ form a projective system which
stabilizes at some $k=k(\lambda)$. In particular, 
$N(\lambda)$ is a finite dimensional space
which coincides with $(N/I^kN)(\lambda)$ for large $k$. 
Thus, the spectrum of ${\bold h}$ in $N$ is bounded from below
(by the minimum of real parts of its eigenvalues on $N/IN$), and
$N=\widehat M$, where $M$ is the direct sum of $N(\lambda)$. 
It is clear that $M\in {\mathcal O}_{c}(\h,G)$. 
The proposition is proved.  
\end{proof}

\begin{remark} Note that Proposition \ref{o=fin} is false for uncompleted
algebras $H_{1,c}(\h,G)$. Proposition \ref{o=fin} is analogous 
to the statement that any 
$D$-module on a formal polydisk $X$ which is finitely generated as 
a module over ${\mathcal O}(X)$ is a multiple of ${\mathcal O}(X)$,
which also fails in the uncompleted case, e.g. for $X={\Bbb A^1}$. 
\end{remark} 

One can also define the completion of $\widehat H_{t,c}(\h,G)_z$ at any point $z\in \h/W$, so that $\widehat H_{t,c}(\h,G)=\widehat H_{t,c}(\h,G)_0$. 
This gives a certain centralizer algebra, as discussed in \cite{BE}, Subsection 3.3. 
See also Subsection \ref{comple1} below for a discussion of this construction in the more general nonlinear case.

\subsection{The Cherednik algebra of a variety with a finite
group action}
 
We will now generalize the definition of $H_{t,c}(\h,G)$
to the case of a nonlinear action of $G$.  
Let $X$ be an affine 
algebraic variety over 
$\Bbb C$, and $G$ be a finite group of automorphisms of $X$.
Let $X^g\subset X$ be the set of fixed points of $g\in G$. 
A component $Y$ of $X^g$ of codimension $1$ in $X$ will be called {\it a reflection hypersurface}.
For instance, if $X=\h$ is a vector space and $G$ acts linearly, then $X^g\subset X$ is a subspace, 
and if ${\rm codim}X^g=1$ then $g$ is a complex reflection and $X^g$ is the corresponding 
reflection hyperplane, which justifies the terminology.   
 
Let $E$ be a $G$-invariant subspace of the space of closed
2-forms on $X$, which projects isomorphically to $H^2(X)$ (it is clear that $E$ exists). 
Consider the algebra \linebreak $G\ltimes {\mathcal O}(T^*X)$, 
where $T^*X$ is the
cotangent bundle of $X$. We are going to define a
deformation $H_{t,c,\omega}(X,G)$ of this algebra
parametrized 
by 

1) complex numbers $t$, 

2) $G$-invariant functions $c$ 
on the (finite) set $S$ of pairs $s=(Y,g)$,
where $g\in G$, $g\ne 1$, and $Y\subset X^g$ is a reflection hypersurface, 
and 

3) elements $\omega\in E^G=H^2(X)^G$. 

If all the parameters are zero, this algebra will coincide with 
$G\ltimes {\mathcal O}(T^*X)$.

Let $t,c=\lbrace{c(Y,g)\rbrace},\omega\in E^G$ be variables. 
Let $D_{\omega/t}(X)_r$ be the algebra (over $\Bbb C[t,t^{-1},\omega]$) of 
twisted (by $\omega/t$) differential operators on $X$ with
rational coefficients. 

\begin{definition}\label{doo}
A Dunkl-Opdam operator for $X,G$ is 
an element of $G\ltimes D_{\omega/t}(X)_r[c]$ given by the formula 
\begin{equation}\label{Dunkl}
D:=t{\bold L}_v+\sum_{(Y,g)\in S}\frac{2c(Y,g)}{1-\lambda_{Y,g}}
\cdot f_Y(x)\cdot (g-1),
\end{equation}
where $\lambda_{Y,g}$ is the eigenvalue of $g$ on 
the conormal bundle to $Y$, $v\in \Gamma(X,TX)$ is a vector field on $X$,
and $f_Y\in \Gamma(X,{\mathcal O}_X(Z))$ 
is an element of the coset $\xi_Y(v)\in \Gamma(X,{\mathcal
O}_X(Z)/{\mathcal O}_X)$ (recall that $\xi_Y$ is defined in
Subsection \ref{agp}). 
\end{definition} 

\begin{definition}
The algebra $H_{t,c,\omega}(X,G)$ is the subalgebra of 
$D_{\omega/t}(X)_r[c]$ generated (over $\Bbb C[t,c,\omega]$) 
by the function algebra 
${\mathcal O}_X$, the group $G$, and the Dunkl-Opdam 
operators. 
\end{definition} 

By specializing $t,c,\omega$ to numerical values, we can define 
a family of algebras over $\Bbb C$, which we will also denote
$H_{t,c,\omega}(X,G)$. Note that when we set $t=0$, the term $t{\bold L}_v$
does not become $0$ but turns into the classical momentum, as in \cite{EM1}, Subsection 
2.10. 

\begin{definition}
$H_{t,c,\omega}(X,G)$ is called the Cherednik algebra of $X, G$. 
\end{definition}

\begin{example}
$X=\h$ is a vector space and $G$ is a 
subgroup in $GL(\h)$. Let $v$ be a constant vector field, 
and let $f_Y(x)=(\alpha_Y,v)/\alpha_Y(x)$, where 
$\alpha_Y\in \h^*$ is a nonzero functional vanishing on $Y$.
Then the operator $D$ is 
just the usual Dunkl-Opdam operator $D_v$ in the complex
reflection case (see Subsection \ref{rca}). 
This implies that all the Dunkl-Opdam operators
in the sense of Definition \ref{doo} 
have the form $\sum f_iD_{y_i}+a$, where $f_i\in \Bbb C[\h]$, $a\in
G\ltimes \Bbb C[\h]$, and $D_{y_i}$ are the usual
Dunkl-Opdam operators (for some basis $y_i$ of $\h$). 
So the algebra $H_{t,c}(\h,G)=H_{t,c,0}(X,G)$ 
is the rational Cherednik algebra for $\h,G$, 
see Subsection \ref{rca}.
\end{example}

The algebra $H_{t,c,\omega}(X,G)$ 
has a filtration $F^\bullet$ which is defined on
generators by $\deg({\mathcal O}_X)=\deg(G)=0$, $\deg(D)=1$ for Dunkl-Opdam
operators $D$. 

\begin{proposition}\label{loca}
Let $\delta\ne 0$ be a $G$-invariant regular function on $X$ which vanishes on reflection hypersurfaces, and let $X^\circ$ be the complement in $X$ of the zero set of $\delta$.  Then we have natural filtration-compatible isomorphisms 
$$
H_{t,c,\omega}(X,G)[\delta^{-1}]\cong H_{t,c,\omega}(X^\circ,G)\cong H_{t,0,\omega}(X^\circ,G)\cong H_{t,0,\omega}(X,G)[\delta^{-1}].
$$
In particular, if $t\ne 0$ then $H_{t,c,\omega}(X,G)[\delta^{-1}]\cong G\ltimes D_{\omega/t}(X^\circ)$ and 
$H_{0,c,\omega}(X,G)[\delta^{-1}]\cong G\ltimes {\mathcal O}(T^*_\omega X^\circ)$. 
\end{proposition} 

\begin{proof} This follows immediately from the definition.  
\end{proof}

\subsection{The formal completion of $H_{t,c,\omega}(X,G)$ at $z\in X/G$}\label{comple1}  

Let us now define the formal completion of $H_{t,c,\omega}(X,G)$ at a point $z\in X/G$, generalizing 
the definition of Subsection \ref{comple} to the nonlinear case. We will assume that the parameters $t,c,\omega$ are numerical; 
the case when they are variables is similar. 

By the Hilbert-Noether theorem, $X/G={\rm Spec}{\mathcal O} (X)^G$ is an affine variety 
(i.e., the algebra ${\mathcal O} (X)^G$ is finitely generated), and $F^iH_{t,c,\omega}(X,G)$ is a finitely generated 
${\mathcal O} (X)^G$-module for each $i$ (under left or right multiplication). 

Let $z\in X/G$. Then we can consider the degree-wise completion 
$$
\widehat{H}_{t,c,\omega}(X,G)_z:=\widehat{{\mathcal O}(X/G)}_z\otimes_{{\mathcal O}(X/G)}H_{t,c,\omega}(X,G), 
$$
where $\widehat{{\mathcal O}(X/G)}_z$ is the algebra of regular functions on the formal neighborhood of $z$ in $X/G$. This is naturally an algebra, with a filtration $F^\bullet$ such that 
$$
F^i\widehat{H}_{t,c,\omega}(X,G)_z=\widehat{{\mathcal O}(X/G)}_z\otimes_{{\mathcal O}(X/G)}F^iH_{t,c,\omega}(X,G).  
$$
Moreover, it turns out that this algebra can be described explicitly via ordinary rational Cherednik algebras. 

Namely, let $x\in X$ be a preimage of $z$, let $G_x\subset G$ be the stabilizer of $x$,
 and let $U_x$ be the formal neighborhood of $x$ in $X$. Then $G_x$ acts on $U_x$. 
By Lemma \ref{cartlem}, any action of a finite group on a formal polydisk over $\Bbb C$ is equivalent to a linear action, 
thus the action of $G_x$ on $U_x$ is equivalent to the linear representation of $G_x$ on $\h:=T_xU_x$. Let $P={\rm Fun}_{G_x}(G,\widehat H_{t,c}(\h,G_x))$ be the space of functions invariant under left multiplication by elements of $G_x$. Then $P$ is a free right $\widehat H_{t,c}(\h,G_x)$-module of rank $N:=|G/G_x|$. 
Let $Z(G,G_x,\widehat H_{t,c}(\h,G_x)):=\End_{\widehat H_{t,c}(\h,G_x)}(P)$ be the corresponding centralizer algebra, see \cite{BE}, Subsection 3.2; it is non-canonically isomorphic to ${\rm Mat}_N(\widehat H_{t,c}(\h,G_x))$. 

Now note that the terms in the Dunkl-Opdam operators \eqref{Dunkl} corresponding to elements 
$g\in G\setminus G_x$ are regular at $x$. Therefore, analogously to \cite{BE}, Theorem 3.2, we have 

\begin{proposition}\label{natisom} 
There is a natural filtered isomorphism 
\begin{equation}\label{filiso}
\widehat{H}_{t,c,\omega}(X,G)_z\cong Z(G,G_x,\widehat H_{t,c}(\h,G_x)). 
\end{equation} 
\end{proposition} 

\subsection{The PBW theorem for $H_{t,c,\omega}(X,G)$}

We will now prove the PBW theorem for $H_{t,c,\omega}(X,G)$. We will assume $t,c,\omega$ are variables. 

\begin{theorem} \label{pbwth} (the PBW theorem) 
We have a natural isomorphism $${\rm gr}_F(H_{t,c,\omega}(X,G))\cong G\ltimes {\mathcal
O}(T^*X)[t,c,\omega].$$ 
\end{theorem}

This implies that for numerical $t,c,\omega$, we have a natural isomorphism $${\rm gr}_F(H_{t,c,\omega}(X,G))\cong G\ltimes {\mathcal
O}(T^*X).$$ 

\begin{proof} 
Suppose first that $X=\h$ is a vector space and $G$ is a 
subgroup in $GL(\h)$. Then, as we mentioned,
$H_{t,c,\omega}(X,G)$ is the rational Cherednik algebra for $G$. 
So in this case the theorem is true, see \cite{EM1}, Subsection 3.2. 

Now consider arbitrary $X$. Let us define a homomorphism of graded algebras
$$
\psi: {\rm gr}_F(H_{t,c,\omega}(X,G))\to G\ltimes {\mathcal O}(T^*X)[t,c,\omega]
$$
({\it the principal symbol homomorphism}). To this end, consider another filtration
$\Phi^\bullet$ on $H_{t,c,\omega}(X,G)$, obtained by restricting 
the usual filtration by order of differential operators from
$G\ltimes D_{\omega/t}(X)_r[c]$ to $H_{t,c,\omega}(X,G)$. Then $F^m\subset \Phi^m$, since the Dunkl-Opdam operators have degree $1$ under $\Phi$. Thus, for any 
$a\in F^m/F^{m-1}$ we can define 
$\psi(a)$ to be the image of $a$ in $\Phi^m/\Phi^{m-1}\subset O(T^*X)_r[t,c,\omega]$ (where $O(T^*X)_r$ are rational functions on $T^*X$ polynomial on the fibers). More explicitly, take $a=gD_{v_1}...D_{v_m}$, $g\in G$, $v_i\in \Gamma(X,TX)$ (any element of $F_m$ is a linear combination of elements of this form). 
Then $\psi(a)=gv_1...v_m$. Hence, $\psi(a)$ in fact belongs to 
$G\ltimes {\mathcal O}(T^*X)[t,c,\omega]$ (i.e., has no poles), i.e., $\psi$ is well defined.    

The homomorphism $\psi$ is surjective, since $G\ltimes O(T^*X)$ is spanned by elements of the form $gv_1...v_m$. Thus, our job is to show that $\psi$ is injective (this is the nontrivial part of the proof). 

Recall that for an affine variety $Y$ and a morphism $f: M\to N$ of finitely generated ${\mathcal O}(Y)$-modules, $f$ is injective iff it is injective on the formal neighborhood of each point of $Y$. Indeed, taking $K={\rm Ker}f$ (a finitely generated ${\mathcal O}(Y)$-module by the Hilbert basis theorem), we find that the completion $K_y$ of $K$ at each $y\in Y$ is zero, which implies that 
the fiber $K|_y$ of $K$ at each $y\in Y$ is zero, giving $K=0$ by Nakayama's lemma. 

In each degree, $\psi$ is a morphism 
of finitely generated modules over ${\mathcal O}(X)^G[t,c,\omega]={\mathcal O}(X/G)[t,c,\omega]$. 
Therefore, to check the injectivity of $\psi$, it suffices to check the injectivity of $\psi$ on the formal 
neighborhood of each point $z\in X/G$. 

Let $x$ be a preimage of $z$ in $X$, and $G_x$ be the stabilizer
of $x$ in $G$. Then $G_x$ acts on the formal neighborhood $U_x$
of $x$ in $X$, and by Lemma \ref{cartlem}, this action is equivalent to the linear action of $G_x$ 
on $\h=T_xU_x$. Therefore, by Proposition \ref{natisom}, the restriction of the map $\psi$ 
to the formal neighborhood of $z$ may be identified with the map 
$$
\widehat\psi: {\rm gr}_FZ(G,G_x,\widehat H_{t,c}(\h,G_x))[\omega]\to Z(G,G_x,G_x\ltimes \widehat {\mathcal O}(T^*\h))[t,c,\omega]
$$
induced by the map 
$$
\overline\psi: {\rm gr}_F\widehat H_{t,c}(\h,G_x)[\omega]\to G_x\ltimes \widehat{\mathcal O}(T^*\h)[t,c,\omega],
$$
where $\widehat{\mathcal{O}}(T^*\h):=\Bbb C[[\h]]\otimes_{\Bbb C[\h]}{\mathcal O}(T^*\h)=\Bbb C[[\h]][\h^*]$. 
Clearly, the injectivity of $\widehat\psi$ is equivalent to the injectivity of $\overline\psi$. 
But the map $\overline\psi$ is obtained by completing at the origin 
the map $\psi$ for $X$ replaced by $\h$ and $G$ by $G_x$.   
Therefore, it suffices to prove the injectivity of $\psi$ in the linear case,
which has been accomplished already. We are done.     
\end{proof}

\begin{remark} The following remark is meant to clarify
the proof of Theorem \ref{pbwth}. 
In the case $X=\h$, the proof of Theorem \ref{pbwth}
is based, essentially, on the (fairly nontrivial) fact 
that the usual Dunkl-Opdam operators $D_v$ commute with each
other. It is therefore very important to note that 
in contrast with the linear case, for a general $X$ we do {\bf
not} have any natural commuting family of Dunkl-Opdam operators. 
Instead, the operators (\ref{Dunkl}) satisfy a weaker property, which is 
still sufficient for the PBW theorem. This property says that 
if $D_1,D_2,D_3$ are Dunkl-Opdam operators corresponding to vector
fields $v_1,v_2, v_3:=[v_1,v_2]$ and some choices of the functions $f_Y$, 
then $[D_1,D_2]-D_3\in G\ltimes {\mathcal O}(X)$ (i.e., it has no
poles). To prove this property, it is sufficient to consider 
the case when $X$ is a 
formal polydisk, with a linear action of $G$. But in this case
everything follows from the commutativity of the ``classical''
Dunkl-Opdam operators $D_v$. 
\end{remark}

\begin{remark} Suppose $G=1$. Then for $t\ne 0$, 
$H_{t,\omega}(X,G)=D_{\omega/t}(X)$. On the other hand,
$H_{0,\omega}(X,G)$ is the Poisson algebra ${\mathcal
O}(T_\omega^*X)$ (in which the Poisson bracket is induced by the $t$-deformation). 
\end{remark} 

\begin{remark} For any $G$, $H_{0,0,\omega}(X,G)=G\ltimes {\mathcal
O}(T^*_\omega X)$, and 
$H_{1,0,\omega}(X,G)=G\ltimes D_\omega(X)$. 
Also, if $\lambda\ne 0$ then $H_{\lambda t,\lambda
c,\lambda \omega}(X,G)=H_{t,c,\omega}(X,G)$. 
\end{remark} 

\begin{remark} The construction of $H_{t,c,\omega}(X,G)$ and 
the PBW theorem extend in a straightforward manner to the case 
when the ground field is not $\Bbb C$ but an algebraically closed
field $k$ of positive characteristic, provided that the order 
of the group $G$ is relatively prime to the characteristic.
\end{remark} 

\subsection{Cherednik algebra as a universal deformation}

\begin{lemma}\label{nont1} 
If for some $(c,\omega)$ the first order deformation 
$H_{1,\hbar c, \hbar\omega}(X,G)$ of \linebreak $G\ltimes D(X)$ 
over $\Bbb C[\hbar]/\hbar^2$ is trivial then $c=0$, $\omega=0$. 
\end{lemma}

\begin{proof} 
Localizing the algebra to formal neighborhoods
of points of $X/G$ using Proposition \ref{natisom}, and using that the lemma is true in the linear
case (Theorem 2.16 in \cite{EG}), we obtain that $c=0$. 
The rest follows from Lemma \ref{nont}. 
\end{proof} 

\begin{theorem}\label{univ}
The algebra $H_{1,c,\omega}(X,G)$ (with formal $c$ and $\omega$)
is a universal formal deformation of $H_{1,0,0}(X,G)=G\ltimes D(X)$.  
\end{theorem}

\begin{proof}
From Corollary \ref{cohGD}, we get: 

\begin{proposition}\label{coho} One has
$$ 
HH^2(G\ltimes D(X),G\ltimes D(X))=H^2(X)^G\oplus(\oplus_{(Y,g)\in
S}H^0(Y))^G.
$$
\end{proposition}

(Note that $H^0(Y)=\Bbb C$; we wrote $H^0(Y)$ to make more
obvious the action of $G$.)

Thus, the dimension of $HH^2(G\ltimes D(X),G\ltimes D(X))$ is the same as the
dimension of the space of parameters $(c,\omega)$. 
Therefore, Theorem \ref{univ} follows immediately from 
Lemma \ref{nont1}. 
\end{proof}

\begin{remark} In fact, it can be shown (by reducing to the linear case using formal completions) 
that the map $$E^G\oplus {\rm Fun}(X,\Bbb C)^G\to HH^2(G\ltimes D(X),G\ltimes D(X))=(\oplus_{(Y,g)\in
S}H^0(Y))^G \oplus H^2(X)^G$$
induced by the deformation $H_{1,c,\omega}(X,G)$ is given by $(c,\omega)\mapsto (c,[\omega])$. 
\end{remark} 

\begin{remark} A special case of the construction of
$H_{t,c,\omega}(X,G)$ is as follows. Let $L$ be a $G$-equivariant
line bundle on $X$. Define $H_{t,c}^L(X,G)$ to be the algebra 
generated by ${\mathcal O}(X)$, $G$, and the Dunkl-Opdam operators
regarded as elements of the smash product of $G$ 
with the algebra of differential operators on $L$
with rational coefficients. It is easy to see that
$H_{t,c}^L(X,G)=H_{t,c,t\omega}(X,G)$, where $\omega$ 
is the curvature of a $G$-stable connection on $L$ (which always exists since $X$ is affine). 
\end{remark} 

\begin{remark} Assume $G=1$. Consider the family of algebras
$H_{t,\omega/t}(X):=H_{t,0,\omega/t}(X,G)$. 
As $t\to 0$, this algebra can be naturally degenerated into 
the algebra $H_\infty(X,\omega)$ defined by 
generators $f\in {\mathcal O}(X)$ and $p_v$, $v\in {\rm Vect}(X)$, 
with defining relations 
$$
fp_v=p_{fv}=p_vf,\ [p_v,p_w]=\omega(v,w).
$$
Thus, $H_\infty(X,\omega)$ is a quantization of the (degenerate)
Poisson structure on $T^*X$, whose Poisson bracket 
is defined by the formula 
$$
\lbrace{f,g\rbrace}=\pi^*\omega(v_f,v_g),
$$
where $f,g\in {\mathcal O}(T^*X)$, and $v_f,v_g$ the
corresponding Hamiltonian vector fields. 
The algebra ${\mathcal O}(T^*X)$ with this Poisson structure 
is obtained as the limit of 
$H_{0,\omega/t}(X)$ when $t$ goes to zero.  
\end{remark}

\begin{remark} The above results generalize without significant
changes to the case when the group $G$ acts on $X$ in a not
necessarily faithful manner. In this case, let $K$ be the kernel
of this action, so that $G/K\subset {\rm Aut}(X)$. Then the
algebra $H_{t,c,\omega}(X,G)$ is defined as above, except that
$$
\omega\in (H^2(X)\otimes \Bbb C[K])^G=(H^2(X)\otimes \Bbb C[K]^K)^{G/K}
$$ 
is a $G/K$-invariant 2-form on $X$ with values in 
the center $\Bbb C[K]^K$ of the group algebra of $K$. The algebra
$H_{t,c,\omega}(X,G)$ satisfies the PBW theorem and 
$H_{1,c,\omega}(X,G)$ is a universal deformation of
$H_{1,0,0}(X,G)=G\ltimes D(X)$. These results are proved similarly
to the case of the faithful action. 
\end{remark} 

\subsection{The spherical subalgebra and the center of the
Cherednik algebra} 

Let $\e=|G|^{-1}\sum_{g\in G}g$ be the symmetrizing idempotent of
$G$. Then we can define the spherical subalgebra 
$\e H_{t,c,\omega}(X,G)\e$ of $H_{t,c,\omega}(X,G)$. 
We denote by $Z_{0,c,\omega}(X,G)$ the center of $H_{0,c,\omega}(X,G)$. 

\begin{theorem} \label{Satake}

(i) ${\rm gr}_F(Z_{0,c,\omega}(X,G))=Z_{0,0,0}(X,G)={\mathcal O}(T^*X/G)$.

(ii) (Satake isomorphism) The map $z\to z\e$ gives an 
isomorphism $Z_{0,c,\omega}(X,G)\to
\e H_{0,c,\omega}(X,G)\e$. In particular, $\e H_{0,c,\omega}(X,G)\e$ is a
commutative algebra. 
\end{theorem}

\begin{proof} (i) First of all, the result holds for $c=0$, since 
$Z_{0,0,\omega}(X,G)={\mathcal O}(T^*_\omega X/G)$. 
In general, the PBW theorem for $H_{t,c,\omega}(X,G)$ implies that we have a natural degree-preserving inclusion 
$\iota: {\rm gr}_F(Z_{0,c,\omega}(X,G))\hookrightarrow Z_{0,0,0}(X,G)$. Assume that $\iota$ is not an isomorphism. 
Then the filtered deformation of $H_{0,0,0}(X,G)$ into $H_{0,c,\omega}(X,G)$ induces a nonzero Poisson bracket $\lbrace,\rbrace$
of some negative degree $d$ on $Z_{0,0,0}(X,G)={\mathcal O}(T^*X/G)$. 
On the other hand, by Proposition \ref{loca} and the special case $c=0$, the corresponding 
inclusion $\iota': {\rm gr}_F(Z_{0,c,\omega}(X^\circ,G))\hookrightarrow Z_{0,0,0}(X^\circ,G)$
is an isomorphism. Hence, the corresponding Poisson bracket $\lbrace{,\rbrace}'$ of degree $d$ on $Z_{0,0,0}(X^\circ,G)={\mathcal O}(T^*X^\circ/G)$ vanishes. But the Poisson bracket $\lbrace,\rbrace$ is the restriction of $\lbrace,\rbrace'$ to a subalgebra, so it must vanish as well. This is a contradiction. Thus, $\iota$ is an isomorphism, proving (i). 

Part (ii) follows from 
(i), since (i) means that the associated graded of the map (ii)
is an isomorphism. 
\end{proof}

Let $M_{c,\omega}(X,G)$ be the spectrum of $Z_{0,c,\omega}(X,G)$. 
It follows from the above that $M_{c,\omega}(X,G)$ is an
irreducible Poisson variety (the Poisson structure comes from the
deformation of $\e H_{0,c,\omega}(X,G)\e$ into $\e H_{t,c,\omega}(X,G)\e$). 
We also have a Lagrangian map $\pi: M_{c,\omega}(X,G)\to X/G$,
whose generic fiber is $T^*_xX$. 

\begin{example} $M_{0,\omega}(X,G)=T^*_\omega X/G$. 
\end{example}

Note that the variety structure of $M_{c,\omega}(X,G)$ is
independent of $\omega$, and only the Poisson bracket depends of
$\omega$. In fact, the dependence of the Poisson bracket on 
$M_{c,\omega}(X,G)$ on $\omega$ is given by the formula 
$$
\lbrace{f,g\rbrace}_\omega=\lbrace{f,g\rbrace}_0+\pi^*\omega(v_f,v_g)
$$
where $v_f,v_g$ are the Hamiltonian vector fields corresponding
to functions $f,g$.

\subsection{Globalization} 
Now let $X$ be any smooth algebraic variety, and $G\subset {\rm Aut}(X)$. 
Assume that $X$ admits a cover by affine $G$-invariant
open sets. Then the quotient variety $X/G$ exists. 

For any affine open set $U$ in $X/G$, let $U'$ be the 
preimage of $U$ in $X$. Then we can define 
the algebra $H_{t,c,0}(U',G)$ as above. 
If $U\subset V$, we have an obvious restriction map
$H_{t,c,0}(V',G)\to H_{t,c,0}(U',G)$. The gluing axiom is clearly satisfied. 
Thus the collection of algebras 
$H_{t,c,0}(U',G)$ can be extended (by sheafification) to a quasicoherent sheaf of algebras on
$X/G$. We are going to denote this sheaf by $H_{t,c,0,X,G}$
and call it the sheaf of Cherednik algebras on $X/G$.
Thus, $H_{t,c,0,X,G}(U)=H_{t,c,0}(U',G)$. 

Similarly, if $\psi\in H^2(X,\Omega_X^{\ge 1})^G$, we can define 
the sheaf of twisted Cherednik algebras $H_{t,c,\psi,X,G}$. 
This is done similarly to the case of twisted differential
operators (which is the case $G=1$). In particular, if 
$L$ is a $G$-equivariant line bundle on $X$, then we can define 
the sheaf of algebras $H_{t,c,X,G}^L$ in the obvious way (it is glued out of
the algebras $H_{t,c}^L(U',G)$).  

Analogously, the varieties $M_{c,\psi}(U',G)$ can be glued into a
single variety $M_{c,\psi}(X,G)$, which is Poisson and has 
a surjective Lagrangian projection $\pi: M_{c,\psi}(X)\to X/G$. 

\begin{definition}
The variety $M_{c,\psi}(X,G)$ is called the Calogero-Moser space of
$X$. 
\end{definition}

{\bf Example.} $M_{0,\psi}(X,G)=T^*_\psi X/G$, the quotient of the twisted 
cotangent bundle of $X$ (\cite{BB}, Section 2) by the $G$-action. 

\begin{remark} These constructions can be generalized to the case when $X$ does not necessarily admit a cover by $G$-invariant affine open sets. 
In this case $X/G$ is an algebraic space which may not be a scheme (\cite{dJ}). But the above constructions go through if instead of $G$-equivariant affine open subsets $U'$ in $X$ (=affine open subsets $U$ in $X/G$) we use $G$-equivariant \'etale morphisms $\pi: U'\to X$, where $U'$ is affine (such a morphism descends to an \'etale morphism $\bar\pi: U\to X/G$, where $U=U'/G$ is affine). 
Thus, we obtain a sheaf of algebras $H_{t,c,\psi,X,G}$ on the algebraic space $X/G$ in the \'etale topology. 
\end{remark} 

\subsection{Modified Cherednik algebras}

Let $X$ is a smooth algebraic variety, and $\psi\in H^2(X,\Omega_X^{\ge 1})$ be a twisting for differential operators on $X$. 
Let $S$ be a finite set of smooth divisors $Y\subset X$, and let $X^\circ$ be the complement of these divisors. 
Let $\eta: S\to \Bbb C$ be a function. 
Denote by $D_{\eta,\psi,X}$ the sheaf of algebras locally generated inside $D_{\psi,X^\circ}$ 
by regular functions and elements $D=\bold L_v+\sum_{Y\in S}\eta(Y)f_Y(x)$, where $v\in \Gamma(X,TX)$ and 
$f_Y$ is a section of $\mathcal{O}_X(Y)/\mathcal{O}_X$ representing $v$. 
Also, let $\psi_Y$ be the class in $H^2(X,\Omega_X^{\ge 1})$ 
defined by the line bundle ${\mathcal O}_X(Y)^{-1}$,
whose sections are functions vanishing
on $Y$. 

\begin{proposition}\label{twi} 
One has a natural isomorophism 
$$
D_{\eta,\psi,X}\cong
D_{\psi+\sum_{Y\in S} \eta(Y)\psi_Y,X}.
$$  
\end{proposition}

\begin{proof} 
The required isomorphism is implemented locally near $x\in X$ by conjugation of twisted differential operators 
by the multivalued function $\prod_{Y\in S} {z_Y^{\eta_Y}}$, where $z_Y$ is a regular function on a neighborhood of $x$ such that the subscheme $Y\subset X$ is cut out in this neighborhood by the equation $z_Y=0$.  
\end{proof}

It will be convenient for us to make a similar modification 
of the sheaf $H_{t,c,\psi,X,G}$. Namely, let 
$\eta$ be a function on the set of conjugacy classes of $Y$ such
that $(Y,g)\in S$. We define $H_{t,c,\eta,\psi,X,G}$ 
in the same way as $H_{t,c,\psi,X,G}$ except that 
the Dunkl-Opdam operators are defined by the 
formula  
\begin{equation}\label{Dunkl1}
D:=t{\bold L}_v+\sum_{(Y,g)\in S}f_Y(x)
\left(\frac{2c(Y,g)}{1-\lambda_{Y,g}}(g-1)+
\eta(Y)\right).
\end{equation}

The following generalization of Proposition \ref{twi} shows that this modification is, in fact,
tautological. 

\begin{proposition} 
One has a natural isomorophism 
$$
H_{t,c,\eta,\psi,X,G}\cong
H_{t,c,\psi+\sum_Y \eta(Y)\psi_Y,X,G}.
$$  
\end{proposition}

\begin{proof} 
Since the terms in the Dunkl-Opdam operators containing $c$ are the same in both cases, it suffices to prove the statement for $c=0$. But this follows from Proposition \ref{twi} (and its classical limit). 
\end{proof}

\begin{remark} We note that the restriction of the sheaf $H_{1,c,\eta,0,X,G}$ to $X^\circ/G$, where $X^\circ:=X\setminus \cup_{(Y,g)\in S} Y$,  is naturally isomorphic to $G\ltimes D_{X^\circ}$. Indeed, this follows from Proposition \ref{loca} and the fact that the 
line bundle ${\mathcal O}_X(Y)$ is trivial on $X^\circ$. 
\end{remark} 

\subsection{Examples} 

In this subsection we will give a few examples of Cherednik
algebras and Calogero-Moser spaces of varieties.
More details on these examples are given in the papers referenced in Section 4. 

\begin{example}\label{curvex} 
Let $C$ be a smooth complex algebraic curve, and 
$X=C^n$, $G=S_n$ ($n>1$). In this case, $c$ is a single parameter. 
If $C=\Bbb C$ then $H_{t,c,0}(X,G)$ is the rational Cherednik
algebra corresponding to the group $S_n$ acting in
the permutation representation $\Bbb C^n$ (see \cite{EG}). 
If $C=\Bbb C^*$ then $H_{t,c,0}(X,G)$ is the trigonometric 
Cherednik algebra for $S_n$ (the degenerate double affine 
Hecke algebra), defined by Cherednik; see e.g. the appendix to \cite{BE}. But for other curves we
obtain new algebras. Note that if the curve is projective, then 
the sheaf $H_{t,c,0,X,G}$ admits a 1-parameter deformation 
$H_{t,c,\psi,X,G}$ (where $\psi$ is a multiple of the class 
of the line bundle $L^{\boxtimes n}$, $L$ being a 
degree $1$ line bundle on $C$). 

Note that for any curve $C$ and $c\ne 0$, the varieties 
$M_{c,\psi}(X,G)$ are smooth and admit a Lagrangian projection to
$S^nC$. This is checked by looking at formal neighborhoods of
points of $X/G=S^nC$, 
where the statement follows from the results of \cite{EG}.

As an example, look at the case $C=\Bbb P^1$, and
 consider the algebra of global sections $A_{t,c,\omega}=
\Gamma(H_{t,c,\omega,X,G})$. It is easy to see that 
$A_{1,0,\psi}=S_n\ltimes D_\psi(\Bbb P^1)^{\otimes n}=
S_n\ltimes {\mathcal D}_\psi^{\otimes n}$, where 
${\mathcal D}_\psi$ is the Dixmier quotient, i.e. the 
quotient of $U(\mathfrak{sl}_2)$ by the central character 
of an $\mathfrak{sl}_2$-module with highest weight $\psi$. 
More generally, let $H_{t,c'}(\Bbb C^n,S_n\ltimes \Bbb Z_2^n)$,
$c'=(c'_1,c'_2)$, be the rational Cherednik algebra of type
$B_n$ (see \cite{BEG1}, Section 6). Let $\bold p$ be the averaging idempotent 
for the group $\Bbb Z_2^n$. Then  
one can show that $A_{t,c,\psi}$ is isomorphic 
to $\bold pH_{t,c'}(\Bbb C^n,S_n\ltimes \Bbb Z_2^n)\bold p$
for $c'$ related to $(c,\psi)$ by the invertible linear
transformation: $c_1'=c, c_2'=\psi-(n-1)c+t/2$. 

This realization allows us to give a Borel-Weil 
type construction of some finite dimensional 
representations of $A_{1,c,\psi}=\bold pH_{1,c'}(\Bbb C^n,S_n\ltimes \Bbb Z_2^n)\bold p$. 
Namely, assume that $t=1$ and $\psi=m$ is a nonnegative integer. 
Then the algebra $A_{1,c,\psi}$ admits a representation $W_m$ of dimension
$(m+1)^n$, on 
global sections of the line bundle $(L^{\boxtimes n})^{\otimes m}$. 
It can be checked that 
this representation has the form $\bold p N(\Bbb C)$, where $N(\Bbb C)$
is a highest weight representation of $H_{1,c'}(\Bbb
C^n,S_n\ltimes \Bbb Z_2^n)$ with trivial highest weight (of
dimension $(2m+1)^n$),
introduced in \cite{BEG1}, Theorem 6.1.   
\end{example}

\begin{example} Let $\h$ be a finite dimensional 
complex vector space, $G$ be a finite subgroup in $PGL(\h)$, and
$G'$ its preimage in $GL(\h)$. Let $\widehat G$ be the intersection
of all subgroups of $G'$ which project onto $G$, and which
contain all complex reflections of $G'$. It is easy to see that
$\widehat G$ is a finite group. Let $K$ be the intersection of $\widehat
G$ with the scalars in $GL(\h)$; 
it is clearly a cyclic group, and $\widehat G/K\subset G$. 
 
Let $X=\Bbb P\h$ be the projective space of $\h$. 
Then we can define the sheaf $H_{t,c,\psi,X,G}$, where 
$\psi$ is a single parameter (corresponding to twisting 
by a power of the tautological line bundle). 

Observe that the set $S$ of pairs $(Y,g)$ can be $G$-equivariantly
identified with the set of complex reflections in $\widehat G$; thus 
the parameter $c$ is
a function on the set of conjugacy classes of complex
reflections in $\widehat G$. 

Consider the algebra of global sections
of the sheaf $H_{t,c,\psi,X,G}$, which we denote 
by $A_{t,c,\psi}$. It is easy to show that 
$A_{1,0,\psi}=G\ltimes D_\psi(\Bbb P\h)$. 
In other words, $A_{1,0,\psi}=(\widehat G\ltimes D(\h)[0])/(E=\psi,K=1)$,
where $D(\h)[0]$ is the algebra of differential operators on $\h$
which commute with the Euler vector field $E$. 
More generally, one can show (see \cite{BM}, Lemma 5.4.1) that 
$$A_{t,c,\psi}=H_{t,c}(\h,\widehat G)[0]/\left({\bold h}=
\frac{\ell t}{2}-\sum_{s\in
S}\frac{2c_s}{1-\lambda_s}+\psi,\ K=1\right),$$ 
(note that ${\bold h}$ is a central element 
of $H_{t,c}(\h,G)[0]$).
\end{example}

\begin{example}
Let us specialize the previous example to the case $\dim \h=2$. 
Thus $G\subset PSL(2,\Bbb C)$. Let $\Gamma$ be the preimage of
$G$ in $SL(2,\Bbb C)$. The kernel of the map from $\Gamma$ to $G$
consists of the identity $1$ and minus identity $z$. 

Note that we have a natural bijection from $SL(2,\Bbb C)\setminus 
\lbrace{1,z\rbrace}$ to the set of pairs $(y,g)$, where $g\in PSL(2,\Bbb C)$ and 
$y\in \Bbb P^1$ is a fixed point of $g$. This map is defined by 
$\gamma\mapsto (y,\bar\gamma)$, where $\bar\gamma$ is the projection
of $\gamma$ to $PSL(2,\Bbb C)$, and $y$ is the fixed line of $\gamma$ on
which it acts with eigenvalue having positive imaginary part. 
This shows that the set $S$ can be identified with 
the set of conjugacy classes in $\Gamma$ of elements not equal to
$1,z$. Let $\bold p_z=(1+z)/2$. 

Consider the sheaf 
$H_{t,c,\psi}(\Bbb P^1,G)$. Let the algebra of its global
sections be denoted by $A_{t,c,\psi}$. It is easy to show 
that $A_{1,0,\psi}=G\ltimes D_\psi(\Bbb P^1)$.

Recall that in \cite{CBH}, Crawley-Boevey and Holland
defined a family of algebras $Q_\lambda$, parametrized by 
elements $\lambda$ of the center of $\Bbb C[\Gamma]$, 
which are the quotient of the smash product 
$\Gamma\ltimes T(\Bbb C^2)$ of $\Gamma$ with the tensor algebra
of the tautological representation by the ideal generated by
$xy-yx-\lambda$ (where $x,y$ is a basis of $\Bbb C^2$).
 
Thus, $A_{1,0,\psi}=\bold p_zQ_\lambda \bold p_z$ for $\lambda=1-z$. 
More generally, one can show that $A_{t,c,\psi}=\bold p_zQ_\lambda \bold p_z$, where
 $\lambda$ is related by an invertible transformation 
with $(t,c,\psi)$ (see \cite{FT}). 

The combination of the above two examples shows that 
the algebra $\bold p_zQ_\lambda \bold p_z$ can be obtained as a quotient of
$H_{t,c}(\h,\widehat G)[0]$.

This allows one to explicitly construct finite dimensional
representations of $\bold p_zQ_\lambda \bold p_z$ (this is done somewhat
implicitly in \cite{CBH}). One way is to 
look at representations of $H_{t,c}(\h,W)[0]$ on 
weight subspaces of representations of $H_{t,c}(\h,W)$ 
from category ${\mathcal O}_{c/t}(\h,W)$ (for $t\ne 0$). 
Another method is to use the isomorphism
$A_{t,c,\psi}\cong\bold p_zQ_\lambda \bold p_z$, and consider finite dimensional
representations of $A_{t,c,\psi}$ for integer $\psi=m$ 
on sections of the $G$-equivariant
vector bundle $L\otimes V$ on $\Bbb P^1$ (where $L$ is the tautological
line bundle, and $V$ is an 
irreducible representation of $\Gamma$ in which $z$ acts by $(-1)^m$). 
\end{example}

\subsection{Finite dimensional $H_{t,c,\psi,X,G}$-modules}

\begin{definition} An $H_{t,c,\psi,X,G}$-module is a quasicoherent
sheaf on $X/G$ with a compatible action of the sheaf of algebras
$H_{t,c,\psi,X,G}$. 
\end{definition} 

Of special interest is the category of finitely
generated $H_{t,c,\psi,X,G}$-modules $M$, i.e.,
quotients of $H_{t,c,\psi,X,G}^{\oplus n}$.
This category is a generalization of the category 
of finitely generated $G$-equivariant twisted D-modules on $X$, 
and deserves a special study, but we will limit
ourselves to discussing finite dimensional objects $M$ in this
category, i.e., such that the space of sections 
$\Gamma(U,M)$ on every open set $U$ is finite dimensional. 
Such modules exist for $t=0$ and also for $t\ne 0$ and special
values of $c$. 

For $x\in X$, let $G_x$ be the stabilizer of $x$ in $G$;
then $G_x\subset GL(T_xX)$. Let 
$c_x$ be the function on the set of conjugacy classes 
of complex reflections in $G_x$ defined by $c_x(g)=c(Y_g,g)$,
where $Y_g$ is the component of $X^g$ that passes through $x$. 

\begin{proposition}
(i) Let $M$ be an indecomposable finite dimensional 
module over $H_{t,c,\psi,X,G}$. Then the set-theoretical support of
$M$ on $X/G$ consists of one point $z\in X/G$. 

(ii) Let $x$ be any point of $X$ projecting to $z$. 
The category of finite dimensional $H_{t,c,\psi,X,G}$-modules
supported at $z$ is equivalent to the category of 
finite dimensional modules over the rational Cherednik algebra
$H_{t,c_x}(T_xX,G_x)$ attached to the group $G_x$ acting in $T_xX$.  
\end{proposition}

\begin{proof}
(i) The statement follows from the fact that the adjoint action 
of ${\mathcal O}_{X/G}(U)$ on $H_{t,c,\psi,X,G}(U')$ is locally
nilpotent. 

(ii) Let $M$ be a finite dimensional $H_{t,c,\psi,X,G}$ module
supported at $z$. Then the 
maximal ideal sheaf $I_z$ of $z$ acts nilpotently on $M$, so 
$M$ can be extended to a module over the degree-wise completion 
$\widehat H_{t,c,\psi,X,G,z}$ of $H_{t,c,\psi,X,G}$
with respect to $I_z$ (see Subsection \ref{comple}). 

Let $\phi: X\to X/G$ be the natural map. 
By Proposition \ref{natisom}, the algebra $\widehat H_{t,c,\psi,X,G,z}$ is isomorphic
(non-canonically) to $\oplus_{y\in \phi^{-1}(z)}\widehat H_{t,c_y}(T_yX,G_y)$,
where hat denotes the completion defined in Subsection \ref{rca}. 
Thus the fiber $M_x$ of $M$ over $x$ (as an
${\mathcal O}$-module) has a natural structure 
of a $\widehat H_{t,c_x}(T_xX,G_x)$-module. 
 Restricting this module to the algebra
$H_{t,c}(T_xX,G_x)$, we obtain a functor 
in one direction, given by $M\mapsto M_x$. 

To construct the functor
in the opposite direction, let $N$ be a finite dimensional module
over $H_{t,c_x}(T_xX,G_x)$. Let $N'={\rm Ind}_{G_x}^G(N)$. 
Arguing as above, we can turn $N'$ into a module over
$H_{t,c,\psi,X,G}$. So the desired functor is $N\mapsto N'$. 
It is easy to verify that the two functors are mutually inverse,
which proves (ii).  
\end{proof}

Thus, the problem of describing finite dimensional representations 
of $H_{t,c,\psi,X,G}$ reduces to the linear case, which is treated  
in \cite{BEG1} and references therein.

\subsection{Quasiinvariants}

In this subsection we will discuss the global version 
of the theory of quasiinvariants for reflection groups. 
One can generalize it to the complex reflection case 
using the ideas of \cite{BC}.

Let $X$ be a smooth complex algebraic variety, $G\subset {\rm
Aut}(X)$, and assume that the quotient $X/G$ exists as a variety. 
Let $S_2\subset S$ be the set of pairs $(Y,g)\subset S$ 
for which $g^2=1$. Any function on $S_2$ can be regarded as a
function on $S$ using extension by zero. 

Recall that $\phi: X\to X/G$ denotes the natural map. 

Let $m: S_2\to \Bbb Z_+$ be a conjugation invariant function. 
The sheaf of $m$-quasiinvariants of $X,G$, 
denoted $Q_{m,X,G}$, is defined to be the
subsheaf of $\phi_*{\mathcal O}_X$ whose local sections satisfy the
condition: for any $(Y,g)\in S_2$, 
$f-{}^gf$ vanishes to order $2m(Y,g)+1$ at $Y$. 
This is a sheaf of rings on $X/G$. This sheaf is the structure
sheaf of an 
algebraic variety $X_m$, together with a natural map 
$\zeta: X\to X_m$ (the normalization map). 
By looking at formal neighborhoods of points, 
it is easy to check (see \cite{BEG})
that $\zeta$ is birational and bijective. 

With these definitions, it is not hard to see that the main
results of \cite{EG1,BEG} extend to the global case. 
Namely, we have the following result. 

\begin{proposition} $X_m$ is a Gorenstein variety. 
\end{proposition}

\begin{proof} The statement is local with respect to $X/G$, so it can be checked 
on formal neighborhoods of points, which is the Feigin-Veselov
conjecture proved in \cite{EG1,BEG}. 
\end{proof}

Now observe that the sheaf of algebras $H_{1,m,0,X,G}$ acts on
$\phi_*{\mathcal O}_X$. Hence the sheaf of spherical subalgebras 
$\e H_{1,m,0,X,G}\e$ acts in $\e\phi_*{\mathcal O}_X=
{\mathcal O}_{X/G}$ by differential
operators. Thus we can extend this action to an action on $\Bbb
C(X)$ (by using the same differential operators). 

\begin{proposition}
$Q_{m,X,G}$ is invariant under the action of the 
sheaf of spherical subalgebras 
$\e H_{1,m,0,X,G}\e$. 
\end{proposition}

\begin{proof}
This is again a local statement, so it can be checked on formal
neighborhoods of points, which is done in \cite{BEG}. 
\end{proof}
 
\section{Cherednik algebras of analytic varieties and Hecke algebras}

The construction and main properties of the Cherednik algebras of
algebraic varieties can be extended 
without significant changes to the case when $X$ is a complex analytic 
variety. We will not specify the 
routine modifications involved; rather, we will use Cherednik
algebras for analytic varieties to define certain  
deformations of orbifold fundamental groups which we call Hecke
algebras, and prove that in certain cases they are flat. 

\subsection{Orbifold fundamental group}\label{s31}
Let $X$ be a connected complex analytic variety, and $G$ is a finite group
of automorphisms of $X$.  
Then $X/G$ is a complex orbifold (for basics on orbifolds, see \cite{Da}). 
Let $x\in X$ be a point with trivial stabilizer. 
In this case we can define the {\it orbifold fundamental group}
$\pi_1^{\rm orb}(X/G,x)$, \cite{Da}. This group is generated 
by homotopy classes of paths on $X$ connecting $x$ and $gx$ for
$g\in G$, with multiplication defined by the rule: $\gamma_1\circ
\gamma_2$ is $\gamma_2$ followed by $g\gamma_1$, where $g$ is
such that $gx$ is the endpoint of $\gamma_2$. 
We have an exact sequence 
$$
1\to \pi_1(X,x)\to \pi_1^{\rm orb}(X/G,x)\to G\to 1.
$$

Note also that $\pi_1^{\rm orb}(X/G,x)=\pi_1((X\times EG)/G,x)$,
where $EG$ is the universal cover of the classifying space of
$G$. In other words, 
the action of $G$ on $X$ gives rise to an associated 
bundle with fiber $X$ over the classifying space $BG$, and 
$\pi_1^{\rm orb}(X/G,x)$ is the fundamental group of 
the total space of this bundle. 

Let $Z$ be the set of all points of $X$ having a nontrivial stabilizer;
then $Z$ is a closed subset of $X$. 
Let $X'=X\setminus Z$.

\begin{definition}  The fundamental group
$\pi_1(X'/G,x)$ will be called the {\it braid group} of
$X/G$. 
\end{definition}

Now let $S$ be the set of pairs $(Y,g)$ such that $g\ne 1$ and 
$Y\subset X^g$ is a reflection hypersurface. 
For $(Y,g)\in S$, let $G_Y$ be the subgroup of $G$ 
whose elements act trivially on $Y$. This group is obviously
cyclic; let $n_Y=|G_Y|$. 
Let $C_Y$ be the conjugacy class in $\pi_1(X'/G,x)$ corresponding
to a small circle going counterclockwise around 
the image of $Y$ in $X/G$. 

\begin{proposition}\label{rel}
The group $\pi_1^{\rm orb}(X/G,x)$ is the quotient 
of the braid group $\pi_1(X'/G,x)$ 
by the relations $T^{n_Y}=1$ for all $T\in C_Y$. 
\end{proposition}

\begin{proof} We have a natural surjective map 
$\theta: \pi_1(X'/G,x)\to \pi_1^{\rm orb}(X/G,x)$ induced by the
embedding $X'\to X$. The kernel of this map obviously 
contains the elements $T^{n_Y}$, $T\in C_Y$.  
 The fact that the kernel is generated by (conjugates of) 
these elements follows from the Seifert-van Kampen theorem.  
\end{proof} 

\subsection{The Hecke algebra of $X,G$}\label{s32}

For any conjugacy class of hypersurfaces 
$Y$ such that $(Y,g)\in S$ we introduce formal parameters
$\tau_{1Y},...,\tau_{n_YY}$. 
The entire collection of these parameters 
will be denoted by $\tau$. 

\begin{definition} We define the Hecke algebra of $X,G$, 
denoted ${\mathcal H}_\tau(X,G,x)$, to be the quotient of 
the $\tau$-adically completed group algebra of the braid group, 
$\Bbb C[\pi_1(X'/G,x)][[\tau]]$, by the $\tau$-adically closed ideal defined by the relations
\begin{equation}\label{poly}
\prod_{j=1}^{n_Y} (T-e^{2\pi ij/n_Y}e^{\tau_{jY}})=0,\ T\in C_Y. 
\end{equation}
\end{definition}

It is clear that up to an isomorphism this
algebra is independent on the choice of $x$, so we will sometimes
drop $x$ form the notation. It follows from Proposition \ref{rel} that 
${\mathcal H}_\tau(X,G)/(\tau=0)=\Bbb C[\pi_1^{\rm
orb}(X/G)]$. Thus, ${\mathcal H}_\tau(X,G)$ is a deformation 
of $\Bbb C[\pi_1^{\rm orb}(X/G)]$. 

\begin{remark} One can also define the Hecke algebra 
${\mathcal H}_\tau(X,G)$ for complex parameters $\tau_{jY}$ (or, rather, $q_{jY}=e^{\tau_{jY}}$),
generalizing the above formal definition in an obvious way.  
\end{remark} 

\subsection{The Knizhnik-Zamolodchikov functor}

In this subsection we will define a global analog of the KZ
functor defined in \cite{GGOR}. 

As we mentioned, similarly to the algebraic case we can define 
the sheaf of algebras $H_{1,c,\eta,0,X,G}$ on $X/G$ (in the analytic topology). 
Note that the restriction of this sheaf to $X'/G$ 
is the same as the restriction of the sheaf $G\ltimes D_X$ to
$X'/G$ (i.e. on $X'/G$, the dependence of the sheaf on the parameters
$c$ and $\eta$ disappears). This follows from the fact 
that the line bundles ${\mathcal O}_X(Y)$ become trivial when
restricted to $X'$. 

Now let $M$ be a module over $H_{1,c,\eta,0,X,G}$
which is a locally free coherent sheaf  
when restricted to $X'/G$. 
Then the restriction of $M$ to $X'/G$ is a
$G$-equivariant D-module on $X'$ which is coherent and locally free 
as an ${\mathcal O}$-module. Thus, $M$ corresponds to 
a locally constant sheaf (local system) on $X'/G$, which 
gives rise to a monodromy representation of the braid group 
$\pi_1(X'/G,x)$ on the fiber $M_x$ of $M$ at $x$. This representation will be denoted by $KZ(M)$. 
This defines a functor KZ, which is analogous to the one in
\cite{GGOR}.

It follows from the theory of D-modules that any ${\mathcal O}_{X/G}$-coherent 
$H_{1,c,\eta,0,X,G}$-module is locally free when restricted to 
$X'/G$. Thus the KZ functor acts from the abelian category 
${\mathcal C}_{c,\eta}$ of ${\mathcal O}_{X/G}$-coherent 
$H_{1,c,\eta,0,X,G}$-modules to the category  
of finite dimensional representations of $\pi_1(X'/G,x)$. It is easy to see that
this functor is exact. 

For any reflection hypersurface $Y$, let $g_Y$ be the generator of $G_Y$ which has
eigenvalue $e^{2\pi i/n_Y}$ in the normal bundle to $Y$. 
Let $(c,\eta)\mapsto \tau(c,\eta)$ be the invertible 
linear transformation defined by the formula
$$
\tau_{jY}=-\frac{2\pi i}{n_Y}\left(2\sum_{m=1}^{n_Y-1}c(Y,g_Y^m)\frac{1-e^{-2\pi
ijm/n_Y}}
{1-e^{-2\pi im/n_Y}}+\eta(Y)\right). 
$$ 

\begin{proposition} \label{Heck} The functor KZ maps
the category ${\mathcal C}_{c,\eta}$ to the category 
of representations of the algebra ${\mathcal
H}_{\tau(c,\eta)}(X,G)$. 
\end{proposition}

\begin{proof} Let $Y$ be a reflection hypersurface. Our job is to show that 
for any $M\in {\mathcal C}_{c,\eta}$, every $T\in C_Y$ satisfies the Hecke relation \eqref{poly} on $KZ(M)$. 

Let $y\in Y$ be a generic point. Using Lemma \ref{cartlem}(i), we may reduce the problem to 
the special case when $X$ is a 1-dimensional open disk $B$ centered at $0$, 
$Y=\lbrace{0\rbrace}$, and $G=G_Y=\Bbb Z/n\Bbb Z$ acts by rotations. 
We may also assume without loss of generality that $\eta=0$, since the algebra
$H_{1,c,\eta}(B,\Bbb Z/n\Bbb Z)$ may be identified with $H_{1,c,0}(B,\Bbb Z/n\Bbb Z)$ by conjugation 
by $z^\eta$, where $z$ is the coordinate on $B$, and this results in rescaling $T$ by 
$e^{2\pi i\eta/n}$. 

Let $M\in {\mathcal C}_{c,0}$, and let $\widehat M:=\Bbb C[[z]]\otimes_{\mathcal O(B)}M$ be the restriction of $M$ 
to the formal neighborhood of zero (where ${\mathcal O}(B)$ is the algebra of holomorphic functions on $B$). Then $\widehat{M}$ is a module over 
$\widehat{H}_{1,c}(\Bbb C,\Bbb Z/n\Bbb Z)$, which is finitely generated over $\Bbb C[[z]]$. By Proposition \ref{o=fin}, 
there exists a unique $M_0\in {\mathcal O}_c(\Bbb C,\Bbb Z/n\Bbb Z)$ 
such that $\widehat M\cong \widehat M_0$, and it is not hard to show 
that $M\cong {\mathcal O}(B)\otimes_{\Bbb C[z]}M_0$ as a $H_{1,c,0}(B,\Bbb Z/n\Bbb Z)$-module
(e.g., for $G=1$ this is just the well known statement that any ${\mathcal{O}}$-coherent $D$-module on $B$ is a multiple of $\mathcal{O}$).
Therefore, we may replace $B$ with $\Bbb C$ and $M$ with $M_0$. 

Now recall \cite{GGOR}, Theorem 5.13:

{\bf Theorem.} If $G\subset GL(\h)$ is a complex reflection group, and 
$N$ is a module over $H_{1,c}(\h,G)$ from category 
${\mathcal O}_c(\h,G)$ then $KZ(N)$ is a module over 
${\mathcal H}_{\tau(c,0)}(\h,G)$. 

Thus, the required statement follows from \cite{GGOR}, Theorem 5.13 for cyclic
groups $G$. 
\end{proof} 

\begin{remark} Note that Theorem 5.13 of \cite{GGOR} for cyclic
groups $G$ is very easy to prove. Namely, let $n=|G|$, and $g$ be the
generator of $G$ whose nontrivial eigenvalue in $\h$ is $e^{2\pi i/n}$. 
Let $\chi_j$ be the character of $G$ given by $\chi_j(g)=e^{2\pi
ij/n}$. An easy computation (see \cite{GGOR}) shows that 
if $N$ is the standard (=Verma) module $M(\chi_j)$ with highest
weight $\chi_j$, then $KZ(N)$ is the 1-dimensional 
character $\zeta_j$ of ${\mathcal H}_\tau$, given by 
$$
\zeta_j(T)=\exp\left(\frac{2\pi
i}{n}\left(j-2\sum_{m=1}^{n-1}c_{g^m}\frac{1-e^{-2\pi
ijm/n}}
{1-e^{-2\pi im/n}}\right)\right).
$$   
This implies the statement for regular $c$, for which the category ${\mathcal O}_c$ is semisimple. 
In general, since $KZ$ is an exact functor, it suffices 
to prove the statement for projective objects $N=P$.
But a projective object $P$ admits a flat deformation to regular values of $c$,
so the statement follows from the regular case by taking a limit.  
\end{remark} 

\subsection{The flatness theorem} 

The main result of this section is the following theorem. 
 
\begin{theorem}\label{p2} Assume that $\pi_2(X)\otimes \Bbb Q=0$. 
Then ${\mathcal H}_\tau(X,G)$ is a {\bf flat} formal deformation 
of $\Bbb C[\pi_1^{\rm orb}(X/G)]$. 
\end{theorem}

The rest of the section contains the proof of Theorem \ref{p2}
and examples of its application. 

\subsection{A lemma on deformations}

We keep the conventions of Subsections \ref{s31}, \ref{s32}. 
Let $\widetilde X$ be the universal covering space of $X$ with base
point $x$. Let $\pi: \widetilde X\to X$ be the covering map. 
Consider the sheaf $\pi_!{\mathcal O}_{\widetilde{X}}$, the direct image with compact
supports of the structure sheaf. Namely, for a small ball $U\subset X$, 
$\Gamma(U,\pi_!{\mathcal O}_{\widetilde{X}})$ is the space of analytic functions 
on $\pi^{-1}(U)$ supported on finitely many connected components of 
$\pi^{-1}(U)$. This sheaf has a natural structure 
of a D-module on $X$ (in general, not coherent). Let $M={\rm Ind}_{D_X}^{G\ltimes D_X}\pi_!{\mathcal O}_{\widetilde{X}}$ be
the corresponding equivariant $D$-module. Thus 
$M$ is a module over $H_{1,0,0,0,X,G}$. 

A central role in the proof of Theorem \ref{p2} is played by 
the following lemma. 

\begin{lemma}\label{defo} If $\pi_2(X)\otimes \Bbb Q=0$ then 
the $G$-equivariant D-module $M$ has a unique flat formal deformation 
to a module over $H_{1,c,\eta,\psi,X,G}$.  
\end{lemma}

\begin{proof} As usual, classes of first order deformations of $M$ lie in
${\rm Ext}^1(M,M)$ and obstructions in ${\rm Ext}^2(M,M)$
(where the Exts are taken in the category of $G$-equivariant $D$-modules on
$X$). So it suffices to show that these two Ext groups vanish. 

Using Shapiro's lemma and the fact that the functor $\pi^*$
(sheaf-theoretic inverse image) is right adjoint to $\pi_!$, we have 
$$
{\rm Ext}^i_{G\ltimes D_X}(M,M)=
{\rm Ext}^i_{D_X}(\pi_!{\mathcal O}_{\widetilde{X}},{\rm Res}M)=
{\rm Ext}^i_{D_{\widetilde X}}({\mathcal O}_{\widetilde{X}},\pi^*{\rm Res}M).
$$
(Here ${\rm Res}$ denotes the functor of forgetting the $G$-equivariant
structure of a $G$-equivariant D-module on $X$).  
But it is clear that $\pi^*{\rm Res}M={\mathcal O}_{\widetilde{X}}\otimes \Bbb
C[\pi_1^{\rm orb}(X/G)]$. Thus, 
$$
{\rm Ext}^i_{G\ltimes D_X}(M,M)=
{\rm Ext}^i_{D_{\widetilde X}}({\mathcal O}_{\widetilde{X}},{\mathcal O}_{\widetilde{X}})\otimes \Bbb
C[\pi_1^{\rm orb}(X/G)].
$$
Thus, by Remark \ref{analvar}, we finally obtain 
$$
{\rm Ext}^i_{G\ltimes D_X}(M,M)=
H^i(\widetilde X,\Bbb C)
 \otimes \Bbb
C[\pi_1^{\rm orb}(X/G)].
$$
Since $\widetilde X$ is simply connected, this clearly vanishes for 
$i=1$. Now consider the case $i=2$. Since $\pi_2(X)\otimes \Bbb Q=0$, 
we have $\pi_2(\widetilde X)\otimes \Bbb Q=0$ (as $\pi_2(X)=\pi_2(\widetilde{X})$), hence by Hurewicz's theorem 
$H^2(\widetilde X,\Bbb C)=0$. Thus ${\rm Ext}^2$ also vanishes, and
we are done.
\end{proof}

\subsection{Proof of Theorem \ref{p2}}

Let $M_{c,\eta,\psi}$ be the flat formal deformation of $M$
whose existence and uniqueness is claimed in Lemma \ref{defo}. 
Then $M_{c,\eta}:=M_{c,\eta,0}$ becomes an ordinary 
$G$-equivariant D-module when restricted 
to the open set $X'$. 

We begin with explaining that even though 
the module $M_{c,\eta}$ 
is, in general, infinitely generated as an
${\mathcal O}$-module (as the cover $\pi$ may have infinitely many sheets), we can 
still apply the functor $KZ$ to it and obtain 
a braid group representation $KZ(M_{c,\eta})$.
Essentially, this is possible because the 
parameters $c,\eta$ are formal, and thus the differential
equation whose monodromy needs to be computed can be solved using
Chen integrals, like in 
Drinfeld's work on the formal KZ equation \cite{Dr}. 

In more detail, 
let $\lbrace{B_i\rbrace}$ be a cover of $X'/G$ by small 
balls. On each $B_i$, the module $M$ can be trivialized, i.e.,
identified with ${\mathcal O}\otimes F$, where $F= \Bbb C[\pi_1^{\rm
orb}(X/G)]$, so that the flat connection on $M|_{X'/G}$ 
becomes the trivial connection. 
Then we get transition maps $g_{ij}^0$ on $B_i\cap B_j$, which are
just elements of $\pi_1^{\rm orb}(X/G)$. 

Now consider the module $M_{c,\eta}$, which is a
deformation of $M_{0,0}=M$. On $X'/G$, this is a deformation 
of (possibly infinite dimensional) bundles with a flat connection, so it can be
understood as a collection of elements $g_{ij}(c,\eta)(z)$
(transition functions), and 
$\omega_i(c,\eta)$ (connection forms).
Here $g_{ij}(c,\eta)$ are formal series in $c,\eta$ 
with coefficients in ${\rm Hom}_{\Bbb C}(F,{\mathcal O}(B_i\cap
B_j)\otimes F)$, and $\omega_i(c,\eta)$ 
are formal series in $c,\eta$ 
with coefficients in ${\rm Hom}_{\Bbb C}(F,\Omega^{1,cl}(B_i)\otimes
F)$, such that $g_{ij}(c,\eta)$ satisfy the cocycle condition for
transition functions, the connections $d+\omega_i(c,\eta)$ 
on $B_i$ glue into a flat connection on $X'/G$, and 
$g_{ij}(0,0)=g_{ij}^0,\omega_i(0,0)=0$. Now, given a path
$\gamma$ starting and ending at $x$, it can be subdivided into
finitely many consecutive segments $\gamma_k$, $k=1,...,n$, each contained in 
a single ball $B_{i_k}$. Here $x$ is the beginning of
$\gamma_1$ and the end of $\gamma_n$, and for convenience we use the same 
$i_1=i_n$ for all $\gamma$. Solving the differential equation 
$df+\omega_{i_k}(c,\eta)f=0$ in $B_{i_k}$ (which we can do using Chen
integrals since $\omega(0,0)=0$, even though ${\rm dim}F$ may be infinite),  
we find the monodromy operator $A_k$ from the beginning to the
end of $\gamma_k$. Clearly, $A_k\in 1+I{\rm End}F$, where $I$ is
the maximal ideal in $\Bbb C[[c,\eta]]$. The monodromy operator 
$A(\gamma)$ along $\gamma$ is then defined to be $A_\gamma=
A_ng_{i_n,i_{n-1}}A_{n-1}...A_2g_{i_2i_1}A_1$. Then
$\gamma\mapsto A_\gamma$ defines a representation of the braid group 
on $F[[c,\eta]]$. This is the desired monodromy representation
$KZ(M_{c,\eta})$.  

Next, we claim that $KZ(M_{c,\eta})$
is, in fact, a representation of the Hecke algebra 
${\mathcal H}_{\tau(c,\eta)}(X,G)$; 
in other words, we can apply Proposition \ref{Heck}
even though the representation may be infinite dimensional. 

To see this, fix a reflection hypersurface $Y$ and a generic point $y\in Y$. 
Now use Lemma \ref{cartlem}(i) to replace $X$ with a
small ball $B$ around $y$ (so that $B$ is invariant under the 
stabilizer $G_Y$ and transversal to $Y$), and $G$ with $G_Y$. We have 
a (non-canonical) isomorphism $M|_B = W \otimes M^B$, where $M^B$ is the analog of the module $M$ for $X$ replaced
with $B$ and $G$ replaced with $G_Y$, and $W$ is a vector space. Namely, $W$ may be identified with the space 
$\Bbb C[\pi_1^{\rm orb}(X/G)/G_Y]$ of finitely supported functions on the homogeneous space $\pi_1^{\rm orb}(X/G)/G_Y$
(here we assume that the base point $x$ is contained in $B$). 
By Lemma \ref{defo},
the unique deformation of $M^B$ to 
a module over $H_{1,c_Y,\eta_Y,0,B,G_Y}$ (where $c_Y,\eta_Y$ are appropriate restrictions of $c,\eta$)
is $M^B_{c_Y,\eta_Y}$. This implies that 
$M_{c,\eta}|_B$ is isomorphic
to the ($(c_Y,\eta_Y)$-adically completed) tensor product $W\otimes M^B_{c_Y,\eta_Y}$. 
But $M^B$ has finite dimensional fibers 
(their dimension is the order $n_Y$ of $G_Y$), 
so Proposition \ref{Heck} can be applied to $M^B$.
Thus we get that $KZ(M_{c,\eta}|_B)$ is a representation of 
the Hecke algebra of $B$ (i.e., the monodromy operator around the center of $B$ satisfies the Hecke relation). Since this is valid for all $Y$, we see that 
$KZ(M_{c,\eta})$ is a representation of the Hecke algebra
${\mathcal H}_{\tau(c,\eta)}(X,G)$, as
desired. 

 Now we are ready to finish the proof of Theorem \ref{p2}. 
As we have shown, the braid group 
representation $KZ(M_{c,\eta})$ factors through the Hecke algebra 
${\mathcal H}_{\tau(c,\eta)}(X,G)$.
If $c,\eta=0$, this representation is the regular representation
of the orbifold fundamental group $\pi_1^{\rm orb}(X/G)$. 
Thus, the regular representation of $\pi_1^{\rm orb}(X/G)$ admits a flat deformation to a representation of the Hecke
algebra ${\mathcal H}_{\tau(c,\eta)}(X,G)$, implying that the Hecke algebra is flat. 
The theorem is proved. 

\subsection{Examples}

In this subsection we would like to discuss a few examples 
of Hecke algebras and of application of Theorem \ref{p2}.

\begin{example}
Let $\h$ be a finite dimensional vector space, and 
$W$ be a complex reflection group in $GL(\h)$. 
Then ${\mathcal H}_\tau(\h,W)$ is the Hecke algebra of $W$
studied in \cite{BMR}. It follows from Theorem \ref{p2} that 
this Hecke algebra is flat. This proof of flatness 
is, in fact, essentially the same as the original proof of this result 
given in \cite{BMR} (based on the Dunkl-Opdam-Cherednik operators). 
\end{example}

\begin{example}
Let $T$ be a maximal torus of a simply connected complex
simple Lie group $G$, and $W=W(T)$ be its Weyl group. 
Then ${\mathcal H}_\tau(T,W)$ is the affine
Hecke algebra. This algebra is also flat by Theorem \ref{p2}. In
fact, its flatness is a well known result from
representation theory; our proof of flatness is essentially due
to Cherednik \cite{Ch}.  
\end{example}

\begin{example}
Let $W,T$ be as in the previous example, and 
$Q^\vee$ be the dual root lattice of $G$. 
Let $E$ be an elliptic curve, and $X=E\otimes Q^\vee$ be the
corresponding Looijenga variety. Then $H_\tau(X,W)$ 
is the double affine Hecke algebra of Cherednik
(\cite{Ch}), and it is flat by Theorem \ref{p2}. 
The fact that this algebra is flat was proved by
Cherednik, Sahi, Noumi, Stokman 
(see \cite{Ch},\cite{Sa},\cite{NoSt},\cite{St})
using a different approach (q-deformed Dunkl
operators). 
\end{example}

\begin{example}\label{Lob}
Let $H$ be a simply connected complex Riemann surface
(i.e., Riemann sphere $\Bbb C\Bbb P^1$, the Euclidean plane $\Bbb C$, or the Lobachevsky plane $\Bbb C_+$), 
and $\Gamma$ be a cocompact lattice in ${\rm Aut}(H)$. 
Let $\Sigma=H/\Gamma$. Then $\Sigma$ is a compact complex
Riemann surface. When $\Gamma$ contains elliptic elements (i.e.,
nontrivial elements of finite order), we are going
to regard $\Sigma$ as an orbifold: it has special points $P_i$, 
$i=1,...,m$ with stabilizers $\Bbb Z_{n_i}$. 
Then $\Gamma$ is the orbifold fundamental group of
$\Sigma$.

Let $g$ be the genus of $\Sigma$, and $a_i,b_i, i=1,...,g$, be the
a-cycles and b-cycles of $\Sigma$. Let $c_j$ be the 
counterclockwise loops around $P_j$.
Then $\Gamma$ is generated by $a_l,b_l,c_j$ with relations
$$
c_j^{n_j}=1,\ c_1c_2...c_m=\prod_la_lb_la_l^{-1}b_l^{-1}.
$$
For each $j$, introduce formal parameters $\tau_{kj}$,
$k=1,...,n_j$. 
Define the Hecke algebra ${\mathcal H}_\tau(\Sigma)$ 
of $\Sigma$ to be generated over $\Bbb C[[\tau]]$ by 
the same generators $a_l,b_l,c_j$ with defining relations 
$$
\prod_{k=1}^{n_j}(c_j-e^{2\pi ij/n_j}e^{\tau_{kj}})=0, \
 c_1c_2...c_m=\prod_la_lb_la_l^{-1}b_l^{-1}. 
$$
  Thus ${\mathcal H}_\tau(\Sigma)$ is a deformation 
of $\Bbb C[\Gamma]$.

We claim that this deformation is flat if $H=\Bbb C$  
or $H=\Bbb C_+$. To show this, let $\Gamma'$ be a normal
subgroup of $\Gamma$ of finite index acting freely on $H$. Let
$G=\Gamma/\Gamma'$, and $X=H/\Gamma'$. Then 
${\mathcal H}_\tau(\Sigma)={\mathcal H}_\tau(X,G)$, so the result
follows from Theorem \ref{p2} and the fact that
$\pi_2(X)=\pi_2(H)=0$. 
\end{example}

\begin{remark}
We note that if $H=\Bbb C\Bbb P^1$  
(so that the condition $\pi_2(X)\otimes \Bbb Q=0$ is violated)
and $\Gamma\ne 1$, then 
this deformation is not flat, which shows that the assumption $\pi_2(X)\otimes \Bbb Q$ in Theorem \ref{p2} cannot be removed. 
Indeed, let $\tau=\tau(\hbar)$ 
be a 1-parameter subdeformation of ${\mathcal H}_\tau(\Sigma)$
which is flat. Let us compute the determinant of the product
$c_1...c_m$ in the regular representation of this algebra (which
is finite dimensional if $H$ is the sphere). On the one hand, it is
$1$, as $c_1...c_m$ is a product of commutators. On the other
hand, the eigenvalues of $c_j$ in 
this representation are 
$e^{2\pi ij/n_j}e^{\tau_{kj}}$ with multiplicity $|\Gamma|/n_j$. 
Computing determinants 
as products of eigenvalues, we get a nontrivial equation on
$\tau_{kj}(\hbar)$, which means that the deformation 
${\mathcal H}_\tau$ is not flat.
\end{remark}

\begin{remark} In the case when $H$ is the Euclidean plane 
(so $\Gamma=\Bbb Z_\ell\ltimes \Bbb Z^2$, $\ell=1,2,3,4,6$), 
the algebras ${\mathcal H}_\tau(\Sigma)$ were considered 
in \cite{EOR}, and it was proved (by a different method) that
they are flat\footnote{Note that closely related flat algebras
(called multiplicative preprojective algebras) were 
considered in \cite{CBS}.}. The case $\ell=1$ is trivial.
In the special case when $\ell=2$, 
this was done earlier by Sahi, Noumi, and Stokman
(\cite{Sa},\cite{NoSt,St}); in this case the Hecke algebra is a
generalized Cherednik's double affine Hecke algebra algebra of rank 1 (for the affine root system $C^\vee C_1$) that ``controls"
Askey-Wilson polynomials.  
\end{remark} 

\begin{remark} If $g=0$, then finite dimensional representations 
of ${\mathcal H}_\tau(\Sigma)$ are closely related to solutions
of the multiplicative Deligne-Simpson problem, see \cite{CB}. 
\end{remark} 

\begin{example}
This is a ``multivariate'' version of the previous example. 
Namely, letting $X,G,\Gamma$ be as in the previous example, 
and $n\ge 1$, we consider the variety $X^n$ with the action of
$S_n\ltimes G^n$. Then ${\mathcal H}_\tau(X^n,S_n\ltimes G^n)$ is
a flat deformation of the group algebra $\Bbb C[S_n\ltimes
\Gamma^n]$ when $X$ is Euclidean or hyperbolic. If $n>1$, this algebra has one more essential parameter 
than for $n=1$; thus we can regard it as a 1-parameter
deformation of $S_n\ltimes {\mathcal H}_{\tau'}(X,G)$, where 
$\tau=(\tau',k)$, and $k$ is a scalar parameter
(in fact, there are two additional parameters, but one of them
is redundant). 
These algebras are considered in the paper \cite{EGO}.
\end{example}

\begin{example}
Recall that in \cite{EG} Ginzburg and the author attach to any 
finite subgroup $G$ in $Sp(2n,\Bbb C)$ a family of algebras
called symplectic reflection algebras, 
parametrized by a complex number $t$ and a function $c$
on the set of conjugacy classes of symplectic reflections in 
$G$. Here we are going to define a lattice analog of this
construction, by assigning a family of algebras parametrized by functions on the set of 
conjugacy classes of affine reflections to any finite subgroup
$G\subset Sp(2n,\Bbb Z)$. 

Let ${\bold L}$ be a symplectic lattice of rank $2n$, and $G$ a finite
subgroup of $Sp({\bold L})$. We would like to deform the group algebra of
the group $\Gamma=G\ltimes {\bold L}$. For this purpose we will realize
$\Gamma$ as the orbifold fundamental group and then pass to the
corresponding Hecke algebra. To do so, let $U=\Bbb R\otimes {\bold L}$
be the corresponding symplectic vector space. Let $\omega$ be the
symplectic form on this space. To make things simple, assume that
$\omega$ is a unique, up to scaling, $G$-invariant symplectic form 
on $U$. Pick a $G$-invariant unitary structure on $U$, i.e., 
a positive definite inner product $(,)$ such that the operator $I: U\to U$ 
defined by 
\begin{equation}\label{unitary}
\omega(v,w)=(Iv,w)
\end{equation}
satisfies the equation $I^2=-1$.
This can be done by suitably normalizing any $G$-invariant positive definite inner product on $U$
 (indeed, the operator $I$ defined by \eqref{unitary} is then skew-adjoint under $(,)$, so $I^2$ 
 is self-adjoint, thus $I^2=\lambda\in \Bbb R$, 
 since $\omega(I^2v,w)$ is another skew-symmetric invariant form on $U$; 
 and $\lambda<0$ since $I$ has imaginary eigenvalues). This makes $U$ into a
complex $n$-dimensional vector space (with complex structure defined by $I$), 
on which $G$ acts $\Bbb C$-linearly. Consider the 
compact complex torus $X=U/{\bold L}$. Clearly, $G$ acts
holomorphically on
$X$, 
$\pi_1^{\rm orb}(X/G)=\Gamma$, 
and $\pi_2(X)=0$, so we can define a flat deformation of $\Bbb
C[\Gamma]$, the Hecke algebra ${\mathcal H}_\tau(X,G)$. 

We will now show that the essential parameters 
in this family of algebras are a complex number $q$ and a
function on the set of conjugacy classes of affine reflections in
$G\ltimes {\bold L}$, i.e., elements whose fixed set in $U$ has
real codimension $2$ (or complex codimension $1$); this is analogous to \cite{EG}. 

Indeed, the parameters $\tau$ correspond by an invertible
transformation to $(c,\eta)$. Recall that $c$
is a function on the set $S$ of pairs $(Y,g)$ such that 
$1\ne g\in G$ and $Y\subset X^g$ is a reflection hypersurface. We observe that $S$ is in a natural bijection 
with the set of conjugacy classes of affine reflections 
  in $G\ltimes {\bold L}$. Thus it remains to observe that 
$\eta$ gives only one more essential parameter, since for any
$Y$, the line bundle ${\mathcal O}_X(GY)$
is $G$-equivariant, and thus its first Chern class is a multiple of
$\omega$. 
\end{example}

\begin{example}\label{ka} (D. Kazhdan). This is a multidimensional version
of Example \ref{Lob}. Let $\widetilde X$ be the $n$-dimensional
complex hyperbolic space $\bold H_{\Bbb C}^n$, 
and $\Gamma$ be a finitely generated 
discrete group of motions of $\widetilde X$,
i.e. a discrete subgroup of $PSU(n,1)$. (There are a number 
of interesting groups of this type generated by complex hypebolic
reflections, such as Mostow groups
$\Gamma(p,t)$, \cite{Mo}). 
Let $\Gamma'$ be a normal subgroup of finite index in $\Gamma$
which acts freely on $\widetilde X$ (it exists by Selberg's
lemma), and let 
$X=\widetilde X/\Gamma'$, $G=\Gamma/\Gamma'$. 
Then the Hecke algebra ${\mathcal H}_\tau(X,G)$ is flat by Theorem \ref{p2}
(since $\pi_2(X)=0$).    
\end{example}

\begin{remark} In fact, in the theory of Cherednik
algebras for analytic varieties discussed in this section, it is not important
that the group $G$ is finite, but it only matters that it acts
on $X$ properly discontinuously. Thus, if $X$ is a complex
manifold and $G$ a discrete group acting properly discontinuously 
and holomorphically on $X$, then one can define the Hecke algebra ${\mathcal
H}_\tau(X,G)$ and show (similarly to Theorem \ref{p2})
that this Hecke algebra is flat if $\pi_2(X)\otimes \Bbb 
Q=0$. Moreover, it is easy
to see that ${\mathcal H}_\tau(X,G)$ is isomorphic to ${\mathcal
H}_\tau(\widetilde X,\pi_1^{\rm orb}(X/G))$. 
(This shows, for instance, that the Hecke algebras of 
Example \ref{ka} are independent on the choice of the group
$\Gamma'$). 
\end{remark} 

\section{A review of further results}

The goal of this section is to describe the developments in the theory of Cherednik and Hecke algebras of varieties since the appearance of  \cite{Et1}.  

\subsection{Representation-theoretic results and generalization of the theory of D-modules to Cherednik algebras} 

The methodology developed in \cite{Et1} was used in \cite{BE} to develop the theory of induction and restriction functors for Cherednik algebras, and in particular (in the appendix to \cite{BE}) to determine the reducibility locus of the polynomial representation of 
the trigonometric Cherednik algebra, giving a new proof of a result of Cherednik. Later, in \cite{BGi} Bellamy and Ginzburg generalized restriction functors of \cite{BE} to Cherednik algebras of varieties, and introduced analogs of Kashiwara's V-filtration and specialization for these algebras.  

The representation theory of $H_{1,c,\psi,X,G}$ was studied more systematically in \cite{Wi}. Namely, this paper determines the possible supports of ${\mathcal O}$-coherent $H_{1,c,\psi,X,G}$-modules, and introduces 
the corresponding category $O_{c,\psi}(X,G)$ of modules over this algebra. This category was 
further studied by D. Thompson in \cite{T1}.  In particular, in this paper he showed that the KZ functor for the algebra $H_{1,c,\psi,X,G}$ is essentially surjective onto the category of finite dimensional representations of the associated Hecke algebra, generalizing a result of Losev \cite{L1} in the linear case. This provides an analog of the Riemann-Hilbert correspondence for sheaves of Cherednik algebras. In the next paper \cite{T2} he generalized to $H_{1,c,\psi,X,G}$ the basic setup of the theory of $G$-equivariant twisted $D$-modules (i.e., the case $c=0$), such as GK dimension, singular support, holonomicity, direct and inverse images, etc. This is again a generalization of the work of Losev in the linear case, \cite{L2}, and of results from \cite{BM}.

\subsection{Deformation-theoretic results} 

Let $M$ be an affine symplectic variety with a symplectic action of a 
finite group $G$. For $g\in G$, let $M^g\subset M$ be the fixed set of $g$. 
Let $A_+$ be a $G$-equivariant formal quantization of $M$
over $\Bbb C[[\hbar]]$, and $A=A_+[\hbar^{-1}]$. It is shown in \cite{DE} 
that 
$$
HH^2(A[G])=HH^2(A^G)=(H^2(M)\oplus \oplus_{g,i: {\rm codim}M_i^g=2}H^0(M^g_i))^G,
$$
where $M^g_i$ are the connected components of $M^g$. It is conjectured in 
\cite{DE} that deformations of $A[G]$ and $A^G$ are unobstructed, i.e., there 
exists a universal deformation parametrized by this $HH^2$. If $M=T^*X$, then this 
follows from Theorem \ref{univ}, so this conjecture is a generalization of Theorem \ref{univ}.  
The conjecture was proved by Halbout and Tang for $G=\Bbb Z/2$ (in the 
similar setting of smooth manifolds), \cite{HT}. 

\subsection{Cherednik algebras of curves} 
In \cite{FG} Finkelberg and Ginzburg studied Cherednik algebras $H_{t,c,\psi,X,G}$ attached to $X=C^n$ and $G=S_n$, 
where $C$ is a smooth algebraic curve (Example \ref{curvex}). In particular, they generalized to the case of an arbitrary curve the 
Hamiltonian reduction construction of the spherical subalgebra $\e H_{t,c,\psi,X,G}\e$, which is given in \cite{EG} in the special case 
$C=\Bbb A^1$. They also studied the Hamiltonian reduction functor from the category of mirabolic character D-modules associated to $C$ 
to the category of representations of this spherical Cherednik algebra, generalizing the results of Gan and Ginzburg in the case of $C=\Bbb A^1$. Later McGerty and Nevins proved a derived microlocalization theorem for these algebras at spherical values of the parameter, \cite{MN}. Also, Calogero-Moser spaces attached to this algebra (Example \ref{curvex}) arise as spaces of ideals in the algebra of differential operators on $C$ (when $C$ is affine), see \cite{Be} and references therein.  

\subsection{Cherednik algebras of abelian varieties and elliptic Calogero-Moser systems for complex crystallographic reflection groups}
 The Cherednik algebra attached to an abelian variety $X$ (for example,
$X=E^n$, where $E$ is an elliptic curve) with an action of a finite
group $G$ is studied in detail in \cite{EM2,EFMV}. Namely, the paper
\cite{EM2} studies a commuting family of Dunkl operators attached to
$X,G$, which act not on functions but rather on sections of line
bundles, i.e. they depend on the ``dynamical" parameter characterizing
a line bundle on $X$ (recall that for a general $X,G$ we don't have
such a family at all). These operators are generalizations to the
complex crystallographic reflection group case of elliptic Dunkl
operators defined in 1994 by Buchshtaber, Felder, and Veselov in the
case of Weyl groups, \cite{BFV}. These operators give rise to certain
representations of $H_{1,c,\psi,X,G}$ which are coherent as
${\mathcal O}_{X/G}$-modules, and application of the generalized KZ
functors to these representations yields monodromy representations of
the Hecke algebra attached to $X,G$ (in particular, of Cherednik's
DAHA in the case of Weyl groups). In \cite{EFMV} elliptic Dunkl
operators were used to prove integrability of elliptic Calogero-Moser
systems for complex crystallographic reflection groups, realizing the
idea of \cite{BFV} (for Weyl groups, this had been proved by Cherednik
in 1995). In the rational and trigonometric case, the quantum
Hamiltonians of the Calogero-Moser system are obtained as symmetric
polynomials of the Dunkl operators, but in the elliptic case this does
not work literally, since Dunkl operators act on sections of line
bundles rather than on functions. Still, it is shown that this
approach works if instead of symmetric polynomials one takes classical Calogero-Moser Hamiltonians, 
in which momentum variables are the elliptic Dunkl operators and position variables are the dynamical parameters labelling the line bundle on which they act. This implies that the elliptic Calogero-Moser system can be described as the algebra of global sections 
of the sheaf $\e H_{1,c,\psi,X,G}\e$ for the critical value of the twist $\psi=\psi(c)$. This is similar 
to the description of the quantum Hitchin system in the geometric Langlands theory, 
as the algebra of global twisted differential operators on the moduli stack of principal bundles on a curve 
for a critical value of twisting. Thus, we obtain quantum integrable systems for any complex crystallographic reflection group, which are new outside the Weyl group case.
 
\subsection{Cherednik algebras of projective spaces and projective quadrics} 
A classical result in the theory of D-modules says that the projective space is D-affine, i.e. the functor of global sections defines an equivalence between the category of D-modules on a projective space and the category of modules over the algebra of global differential operators; this is a special case of the Beilinson-Bernstein localization theorem. In other words, the sheaf of differential operators on a projective space is affine. Moreover, this generalizes to the sheaf of twisted differential operators when the twisting parameter is not a negative integer. 
The paper \cite{BM} generalizes this result to Cherednik algebras of the projective space, $H_{1,c,\psi,\Bbb P(V),G}$, where $V$ is a linear representation of a finite group $G$ (so that the classical result is recovered for $c=0$). Namely, it shows that the sheaf $H_{1,c,\psi,\Bbb P(V),G}$ for generic parameters $c,\psi$ is affine, and determines explicitly the exceptional set. It is also interesting to consider the sheaf 
$H_{1,c,\omega,Q,G}$, if $G$ preserves a nondegenerate quadratic form $\bold q$ on $V$, and $Q\subset \Bbb P(V)$ 
is the quadric $\bold q=0$. In this case, it turns out that the algebra of global sections of $H_{1,c,\psi,Q,G}$
coincides with the central reduction of the {\it Dunkl angular momenta algebra} defined in \cite{FH}. This is discussed in much more detail in \cite{FT}.      
 
\subsection{Algebraic PBW theorems for Hecke algebras of complex manifolds with a finite group action}

Let $X$ be a connected complex manifold and $G$ be a finite group of holomorphic transformations of $X$. 
Then Theorem \ref{p2} says that the Hecke algebra ${\mathcal{H}}_\tau(X,G)$ saisfies the formal PBW theorem, provided that
$\pi_2(X)\otimes \Bbb Q=0$. However, the Hecke algebra has an algebraic version ${\mathcal{H}}(X,G,\bold q)$ over $\Bbb C[\bold q,\bold q^{-1}]$, where $\bold q=(q_{jY})$, $q_{jY}=e^{\tau_{jY}}$, and one may ask if in the same situations we have the algebraic PBW theorem: the algebra ${\mathcal{H}}(X,G,\bold q)$ is a free module over $\Bbb C[\bold q,\bold q^{-1}]$. This is a much stronger statement than the formal PBW theorem given by Theorem \ref{p2}, and it is known in very few situations (apart from the classical cases of finite, affine, and double affine Hecke algebras associated to Weyl groups). Namely, for $X$ a vector space and $W$ acting linearly, this is the Brou\'e-Malle-Rouquier conjecture, \cite{BMR}, which is now is finally proved (over $\Bbb C$), see \cite{Et2}. Also, if $X$ is a compact Riemann surface of genus $\ge 1$, the algebraic PBW theorem is proved in \cite{EOR,ER1}. Finally, the paper \cite{ER2} studies the case $X=\Bbb \Bbb C P^1$, in which $\pi_2(X)\otimes \Bbb Q\ne 0$ and the PBW theorem fails (already at the formal level), and determines the flatness locus of the Hecke algebra in the space of parameters.

\end{document}